\documentclass[10pt,a4paper]{amsart} 
\usepackage{amsfonts,amsmath,amssymb,amsthm,arydshln,dsfont,enumitem,ifpdf,ifthen,mathtools,slashed,xargs,xparse}

\usepackage[utf8]{inputenc} 
\usepackage[T1]{fontenc}

\usepackage[english]{babel}
\usepackage{csquotes}

\usepackage[hidelinks]{hyperref}
\usepackage[ddmmyyyy]{datetime}

\usepackage{tikz}
\usetikzlibrary{matrix,arrows}
\usepackage{tikz-cd}

\usepackage[overload]{textcase}

\ifpdf
  \usepackage[final,expansion=alltext,protrusion=alltext]{microtype} 
\fi

\theoremstyle{plain}

\newtheorem*{Th*}{Theorem}

\newtheorem*{Cor*}{Corollary}

\theoremstyle{definition}

\theoremstyle{remark}

\newif\ifTNS 
\TNStrue

\def\printtheoremname#1{\csname#1name\endcsname}
\def\printtheoremnames#1{\csname#1names\endcsname}

\def\thmref#1#2{\printtheoremname{#1}\ifTNS~\fi\ref{#1:#2}}

\def\uc#1#2{\MakeUppercase{#1}{#2}} 

\newcommand{\DefTheorem}[2]{\newenvironmentx{#1}[2][1=\empty,2=\empty]{%
    \ignorespaces%
    \ifx##2\empty%
      \begin{#2}%
    \else%
      \begin{#2}[{\uc##2}]%
    \fi%
    \ifx##1\empty%
      {}%
    \else%
      \label{#1:##1}%
    \fi%
    \ignorespaces}{\end{#2}\ignorespacesafterend}}

\newcommand{\prfof}[2]{\protect{Proof of~\thmref{#1}{#2}}}

\makeatother

\DefTheorem{Th}{theorem}
\DefTheorem{Prop}{proposition}
\DefTheorem{Cor}{corollary}
\DefTheorem{Lem}{lemma}
\DefTheorem{Def}{definition}
\DefTheorem{Rem}{remark}
\DefTheorem{Par}{para}
\DefTheorem{Not}{notation}
\DefTheorem{Exer}{exercise}
\DefTheorem{Ex}{example}
\DefTheorem{Cons}{construction}
\DefTheorem{Scho}{scholium}

\newenvironment{Par*}{\ignorespaces\noindent\ignorespaces}{\ignorespacesafterend}

\numberwithin{equation}{section}

\newcommand\Define[2][\empty]{\ignorespaces%
  \emph{#2}}%

\tikzset{
  commutative diagrams/.cd,
  arrow style=tikz,
  diagrams={>=stealth},
  shift up/.style={
    to path={([yshift=#1]\tikztostart.east) -- ([yshift=#1]\tikztotarget.west) \tikztonodes}},
  shift up left/.style={
    to path={([yshift=#1]\tikztostart.west) -- ([yshift=#1]\tikztotarget.east) \tikztonodes}},
  mathdouble/.style={-,double equal sign distance}
}

\def\ger{\mathfrak}

\newcommand\SheafTypeface{\mathcal}
\def\sh{\SheafTypeface}

\def\DMO{\DeclareMathOperator}

\newcommand\ev{{\bar 0}}
\newcommand\odd{{\bar 1}}

\newcommand{\defi}{\coloneqq}     

\newcommand{\Fa}{For all }
\newcommand{\fa}{for all }

\newcommand{\fs}{for some }

\newcommand{\scth}{such that }

\newcommand\cf{\emph{cf.}~}

\newcommand\ie{\emph{i.e.}~}

\newcommand\opcit{\emph{op.~cit.}}

\newcommand\vphi{\varphi}
\newcommand\vrho{\varrho}

\newcommand\eps{\varepsilon}
\newcommand\nats{\mathbb{N}}
\newcommand\ints{\mathbb{Z}}

\newcommand\reals{\mathbb{R}}
\newcommand\cplxs{\mathbb{C}}

\def\aff{{\mathbb A}}

\newcommand\sle{\leqslant}
\newcommand\sge{\geqslant}

\DMO\dom{\mathrm{dom}}
\DMO\rk{\mathrm{rk}}
\DMO\cork{\mathrm{cork}}
\DMO\Ad{\mathrm{Ad}}
\DMO\ad{\mathrm{ad}}
\DMO\GL{\mathrm{GL}}
\DMO\id{\mathrm{id}}
\DMO\pr{\mathrm{pr}}
\DMO\gr{\mathrm{gr}}
\DMO\sll{\ger{sl}}
\DMO\sdim{\mathrm{sdim}}
\DMO\sgn{\mathrm{sgn}}
\DMO\re{\mathrm{Re}}
\DMO\coker{\mathrm{coker}}
\DMO\im{\mathrm{im}}
\DMO\coim{\mathrm{coim}}
\DMO\codim{\mathrm{codim}}
\DMO\supp{\mathrm{supp}}
\DMO\str{\mathrm{str}}
\DMO\tr{\mathrm{tr}}
\DMO\res{\mathrm{res}}
\DMO\vol{\mathrm{vol}}

\makeatletter
\newcommand\Size[7][1]{
 \ifx#20%
        \def\r@l{}\def\r@m{}\def\r@r{}%
 \else%
    \ifx#21%
           \def\r@l{\bigl}\def\r@r{\bigr}\def\r@m{\bigm}%
    \else%
           \ifx#22%
                 \def\r@l{\Bigl}\def\r@r{\Bigr}\def\r@m{\Bigm}%
            \else%
                 \ifx#23%
                        \def\r@l{\biggl}\def\r@r{\biggr}\def\r@m{\biggm}%
                  \else
                        \ifx#24%
                        \def\r@l{\Biggl}\def\r@r{\Biggr}\def\r@m{\Biggm}%
                        \fi%
                  \fi%
            \fi%
      \fi%
 \fi%
 \ifx#10%
       \def\r@m{}%
 \fi%
 \r@l#3{#4}\r@m#5{#6}\r@r#7%
}%

\makeatother

\makeatletter
\def\Set@Scallop[#1]#2#3{{#1}\Parens{#2}{#3}}
\newcommand\DeclareScalableOperator[2]{%
  \expandafter\def\csname#1\endcsname{\@ifnextchar[{{#2}\Set@Scallop}{{#2}\Set@Scallop[{}]}}
}
\makeatother

\def\DSO{\DeclareScalableOperator}

\DSO{Presh}{\cat{Presh}}
\DSO{Sh}{\cat{Sh}}

\DSO{ShHom}{\sh Hom} 
\DSO{ShGHom}{\underline{\sh Hom}} 
\DSO{ShEnd}{\sh End} 
\DSO{ShGEnd}{\underline{\sh End}} 
\DSO{ShDer}{\sh Der} 
\DSO{ShGDer}{\underline{\sh Der}} 
\DSO{ShCliff}{\sh C\text{\emph{liff}}} 
\DSO{Db}{\sh Db}

\DSO{Ct}{\mathcal{C}} 
\DSO{Hom}{\mathrm{Hom}} 
\DSO{End}{\mathrm{End}} 
\DSO{Aut}{\mathrm{Aut}} 
\DSO{ABer}{\Abs0{\mathrm{Ber}}}
\DSO{Ber}{\mathrm{Ber}}
\DSO{Coind}{\mathrm{Coind}}
\DSO{Der}{\mathrm{Der}}
\DSO{GDer}{\underline{\mathrm{Der}}}
\DSO{GHom}{\underline{\mathrm{Hom}}} 
\DSO{GEnd}{\underline{\mathrm{End}}} 
\DSO{PW}{\mathrm{PW}} 

\newcommand\Set[3]{
                                 \Size{#1}{\{}{#2}{|}{#3}{\}}%
}%
\newcommand\Dual[3]{
                                 \Size[0]{#1}{\langle}{#2}{,}{#3}{\rangle}%
}%
\newcommand\Parens[2]{
  \Size[0]{#1}{(}{#2}{}{}{)}
}
\newcommand\Bracks[2]{
  \Size[0]{#1}{[}{#2}{}{}{]}
}

\newcommand\Norm[2]{
  \Size[0]{#1}{\lVert}{#2}{}{}{\rVert}
}
\newcommand\Abs[2]{
  \Size[0]{#1}{\lvert}{#2}{}{}{\rvert}
}

\makeatletter

\makeatletter

\newif\if@smallmat
\newif\if@none
\newif\if@paren
\newif\if@brack
\newif\if@brace
\newif\if@vline
\newenvironment{Matrix}[2][1]
{\ifx#20%
      \@smallmattrue%
\else%
       \@smallmatfalse
\fi%
\ifx#11%
       \@nonefalse\@parentrue\@brackfalse\@bracefalse\@vlinefalse%
\else%
     \ifx#12%
          \@nonefalse\@parenfalse\@bracktrue\@bracefalse\@vlinefalse%
      \else%
          \ifx#13%
               \@nonefalse\@parenfalse\@brackfalse\@bracetrue\@vlinefalse%
          \else%
               \ifx#14%
                     \@nonefalse\@parenfalse\@brackfalse\@bracefalse\@vlinetrue
               \else%
                     \ifx#15%
                           \@nonefalse\@parenfalse\@brackfalse\@bracefalse\@vlinefalse%
                     \else%
                           \@nonetrue\@parenfalse\@brackfalse\@bracefalse\@vlinefalse%
                     \fi%
               \fi%
          \fi%
      \fi%
 \fi%
 \if@smallmat%
      \if@none%
           \begin{smallmatrix}%
      \else%
          \if@paren%
                \left(\begin{smallmatrix}%
          \else%
                \if@brack%
                        \left[\begin{smallmatrix}%
                \else%
                        \if@brace%
                             \left\{\begin{smallmatrix}%
                        \else%
                             \if@vline%
                                  \left\lvert\begin{smallmatrix}%
                              \else%
                                  \left\lVert\begin{smallmatrix}%
                              \fi%
                        \fi%
                \fi%
          \fi%
      \fi%
 \else%
      \if@none%
           \begin{matrix}%
      \else%
          \if@paren%
                \begin{pmatrix}%
          \else%
                \if@brack%
                        \begin{bmatrix}%
                \else%
                        \if@brace%
                             \begin{Bmatrix}%
                        \else%
                             \if@vline%
                                  \begin{vmatrix}%
                              \else%
                                  \begin{Vmatrix}%
                              \fi%
                        \fi%
                \fi%
          \fi%
      \fi%
 \fi}%
{\if@smallmat%
      \if@none%
           \end{smallmatrix}%
      \else%
          \if@paren%
                \end{smallmatrix}\right)%
          \else%
                \if@brack%
                        \end{smallmatrix}\right]%
                \else%
                        \if@brace%
                             \end{smallmatrix}\right\}%
                        \else%
                             \if@vline%
                                  \end{smallmatrix}\right\rvert%
                              \else%
                                  \end{smallmatrix}\right\rVert%
                              \fi%
                        \fi%
                \fi%
          \fi%
       \fi%
 \else%
      \if@none%
           \end{matrix}%
      \else%
          \if@paren%
                \end{pmatrix}%
          \else%
                \if@brack%
                        \end{bmatrix}%
                \else%
                        \if@brace%
                             \end{Bmatrix}%
                        \else%
                             \if@vline%
                                  \end{vmatrix}%
                              \else%
                                  \end{Vmatrix}%
                              \fi%
                        \fi%
                \fi%
          \fi%
      \fi%
 \fi}%

\makeatother

\def\clap#1{\hbox to 0pt{\hss#1\hss}} 

\begin{document}
  \title{Non-Euclidean Fourier inversion on super-hyperbolic space}

  \author[Alldridge]
  {Alexander Alldridge}

  \address{Universit\"at zu K\"oln\\
  Mathematisches Institut\\
  Weyertal 86-90\\
  50931 K\"oln\\
  Germany}
  \email{alldridg@math.uni-koeln.de}

  \author[Palzer]
  {Wolfgang Palzer}

  \thanks{Research funded by DFG, grant nos.~AL 698/3-1, TRR 12, TRR 183, and ZI 513/2-1, and the Institutional Strategy of the University of Cologne in the Excellence Initiative}

  \subjclass[2010]{Primary 43A30, 43A85; Secondary 22E30, 22E46, 43A90, 58A50}

  \keywords{Lie supergroup, Fourier inversion formula, Harish-Chandra $c$-function, Non-Euclidean Helgason--Fourier transform, spherical function, Riemannian symmetric superspace}

  \begin{abstract}
  	For the super-hyperbolic space in any dimension, we introduce the non-Euclidean Helgason--Fourier transform. We prove an inversion formula exhibiting residue contributions at the poles of the Harish-Chandra $c$-function, signalling discrete parts in the spectrum. The proof is based on a detailed study of the spherical superfunctions, using recursion relations and localization techniques to normalize them precisely, careful estimates of their derivatives, and a rigorous analysis of the boundary terms appearing in the polar coordinate expression of the invariant integral. 
  \end{abstract}

  \maketitle

  \section*{Introduction}

  The super-hyperbolic space $\mathrm{SOUSp}_0(1,1+p\mid 2q)/\mathrm{SOSp}(1+p\mid 2q)$, $p>0$, is a Riemannian symmetric supermanifold generalizing the Riemannian hyperboloid
  \[
  	H_p\defi\mathrm{SO}_0(1,1+p)/\mathrm{SO}(1+p)=\Set1{x\in\reals^{1+p}}{x_0>0,x_0^2-x_1^2-\dotsm-x_p^2=1}.
  \]
  In the framework of S.~Helgason's theory of harmonic analysis on Riemannian symmetric spaces \cite{helgason-gga,helgason-sym}, one considers non-Euclidean generalizations of the classical Fourier transform. In the case of the hyperboloid, it maps functions $f(x)$ on $H_p$ to functions $\sh F(f)(\lambda,b)$ on $\cplxs\times\mathbb S^p$. The inversion formula for $\sh F$ states that 
  \begin{equation}\label{eq:fi-classical}
  	Cf(x)=\sh J(\sh F(f))(x)\defi\int_{-\infty}^\infty\frac{ds}{\Abs0{\mathbf c(is)}^2}\int_Bdb\,\sh F(f)(is,b)\Dual0xb^{-is-\vrho}  	
  \end{equation}
  where $B=\{1\}\times\mathbb S^p$, $\Dual0\cdot\cdot$ is the standard bilinear form of signature $(1,p)$, $C$ is some positive constant, $\vrho=\frac p2$, and $\mathbf c(\lambda)$ is Harish-Chandra's $c$-function. 

  This formula can be given an interpretation in terms of the representation theory of $G=\mathrm{SO}_0(1,1+p)$. The Lie group $G$ acts transitively by isometries on the hyperbolic space $H_p=G/K$ (where $K=\mathrm{SO}(1+p)$). Thus, the regular representation on $L^2(H_p)=L^2(G/K)$ is unitary. The inversion formula expresses this representation as the multiplicity-free direct integral of the unitary spherical principal series representations of $G$, which are parametrised by $\lambda=is\in i\reals$ and realised on sections of line bundles on $B=K/M$.

  The topic of this paper is to investigate to which extent the inversion formula for the non-Euclidean Fourier transform on hyperbolic space generalizes to the super case.
  Supermanifolds such as the super-hyperboloids appear as the target spaces of non-linear $\sigma$-models that have been applied extensively in physics, in the study of localization and delocalization in disordered metals, semiconductors, and superconductors \cite{alzi}, and more recently, in the context of symmetry-protected topological phases of matter \cite{bagrets-altland}. By general procedures, statistics of ensembles of random quantum Hamiltonians can be related to solutions of geometric PDE on the target supermanifold in the $\sigma$-model approximation, and it is in their study that harmonic analysis comes into play. 

  The super-hyperboloid, in particular, occurs as a toy model (ignoring the ``compact sector'') for the $\sigma$-model of class $BD\mathrm I|C\mathrm{I\!I}$, in the parlance of Altland--Zirnbauer \cite{alzi,hhz}. (The corresponding quantum Hamiltonians lie in class $A\mathrm I$.) In the simplest case, where $p=q=1$ and hence $\vrho=\frac p2-q=-\tfrac12$, the super-hyperboloid is just the super-hyperbolic disc studied some time ago by M.~Zirnbauer \cite{zirn-cmp}. As he showed, the harmonic analysis on this space exhibits a number of striking peculiarities; notably, an additional term appears in the inversion formula \eqref{eq:fi-classical}. Moreover, due the fact that $\vrho$ is negative, the Riemannian volume exhibits exponential decay at infinity (instead of exponential growth) and, when expressed in polar coordinates, a singularity at the origin. This analysis leads to the precise prediction of the transition from a diffusive regime to one of exponential localization in a thin wire, as a function of system size in units of correlation length.

  \medskip\noindent
  In this paper, we prove a Fourier inversion formula for the super-hyperboloid, for any choice of $p>0$ and $q\sge0$ (\thmref{Th}{fouinv}). As it turns out, the formula is a function of $\vrho=\frac p2-q$ alone; while for $\vrho\sge 0$, it takes the same form as in the classical case, for $\vrho<0$, there are additional contributions to Equation \eqref{eq:fi-classical}: 
  \begin{equation}
      2^{2(1-\vrho)}\pi f=
      \begin{cases}
        \sh J(\sh F(f)), &\text{if }\vrho\sge0,\\
        \sh J(\sh F(f))-f*\sh J(1), &\text{if }\vrho<0.
      \end{cases}
  \end{equation}
  Here, $*$ denotes the convolution product induced by the action of the Lie supergroup $G=\mathrm{SOUSp}_0(1,1+p|2q)$. In case $\vrho<0$, $\sh J(1)$ exists due to the exponential decay of the volume, and is given by the residues at a finite number of  poles of the Harish-Chandra $c$-function. (The resulting formula is somewhat reminiscent of the case of the non-Riemannian hyperboloids \cite{andersen}.) 

  At present, a precise interpretation of the inversion formula in terms of the representation theory of the Lie supergroup $G$ is not available. One reason is that $G$ \emph{is not a real Lie supergroup}---a real form is only fixed for the underlying Lie group. In fact, this is unavoidable as there is (for generic $p$ and $q$) no real form of the Lie superalgebra $\ger{osp}(2+p|2q,\cplxs)$ whose even part acts by infinitesimal isometries on the super-hyperboloid. Thus, there is no notion of unitary representations and no obvious generalization of the Hilbert space $L^2(H_p)$ at hand. Nonetheless, one can say that the ``most continuous'' part of the inversion formula (\ie~$\sh J(\sh F(f))$) is given by spherical principal series representations of $G$ unitary when restricted to the underlying Lie group $G_0=SO_0(1,1+p)\times\mathrm{USp}(2q)$. The term $f*\sh J(1)$ appearing for $\vrho<0$ is a discrete contribution by a finite number of spherical representations. 

  Let us comment on the proof of our main result. As in the case studied by Zirnbauer, the additional term in the inversion formula is related to the occurrence of ``boundary terms'' in the polar coordinate expression of the Riemannian Berezin integral on the super-hyperboloid. In general, the boundary terms have a more intricate form than in Zirnbauer's case, and we give the general expression in \thmref{Th}{inv-ber-gk}. Notably, when $\vrho$ is negative, the integer and half-integer cases lead to radically different integro-differential expressions.

  The fundamental dichotomy between $\vrho$ integral or half-integral pervades the entire article, and is closely related to a similar one for the hyperbolic spaces $H_p$. When $\vrho$ is a negative integer, the inversion formula can be reduced to a study of the $K$-invariant case (the spherical transform), and, by a recursion on $\vrho$, to the case where $\vrho\sge0$. By contrast, in the half-integral case, such a reduction is not possible, as the boundary supersphere $B=K/M=\mathrm{SOSp}(1+p\mid 2q)/\mathrm{SOSp}(p\mid 2q)$ then has volume zero. Rather, the proof is based on a delicate analysis of boundary terms. 

  The recursive procedure valid for integral $\vrho<0$ can also be applied when $\vrho\sge0$, and we prove the inversion formula in this case without using Helgason's general result, reducing it instead to the inversion formula for the Euclidean Fourier transform (in case $\vrho$ is integral) and the Mehler--Fock transform (in the half-integral case). This allows us to deduce our main result in a pedestrian fashion, virtually without recourse to the general theory. Our approach is closely related to the so-called shift operators introduced by Opdam.

  Throughout the article, we have taken great pains to normalize all quantities with ultimate precision, in particular, the Berezinian densities and the spherical superfunctions. This is a very delicate matter, especially when $B$ has volume zero (namely, when $\vrho<0$ is half-integral), and adds substantially to the length of the exposition. However, we believe this to be a valuable piece of information, since at the outset, it was not all obvious what the correct normalization should be. 

  Let us make some final remarks to the situation for more general Riemannian symmetric superspaces $X$, in particular of high rank. A serious complication is that even when $X$ is irreducible, the underlying Riemannian symmetric space $X_0$ will usually be the product of spaces of the non-compact and the compact type, precluding an easy generalisation of the rank-one theory. As yet unexplored is the relation to a super version of Helgason's Radon transform. It seems probable that such a relation exists, and it may be useful to exploit.

  \medskip\noindent
  This work is based on our previous results \cite{ap-sph} on the $c$-function and the Harish-Chandra expansion of spherical superfunctions on rank-one Riemannian symmetric superspaces of non-compact type. The results on the Harish-Chandra expansion are used in Section \ref{sec:wavepacket} to derive estimates on the spherical superfunctions and their derivatives, and to study the residues of the wave packet transform $\sh J$, in particular, when applied to the constant function $1$. The results on the $c$-function, however, are rederived by more elementary means in Section \ref{sec:recursion}, as a byproduct of our recursion relation for the spherical superfunctions. Following these preliminaries, we prove, in Section \ref{sec:inversion}, the polar coordinate expression of the invariant integral on the super-hyperboloid, and then, the inversion formula. Besides the facts already mentioned, our results from Ref.~\cite{ap-sph} on the localization of supersphere integrals are also applied, in order to evaluate the boundary terms in case $\vrho<0$ is half-integral, and to determine the normalization of Berezinians in this case.


  \medskip\noindent
  \emph{Acknowledgements.} First and foremost, we wish to thank Martin Zirnbauer for many insightful discussions. We also thank two anonymous referees for their comments and suggestions for improvement. Moreover, the first-named author gratefully acknowledges the hospitality of the Institute for Theoretical Physics at the University of Cologne during the preparation of the article. Large parts of the results presented here were obtained in the second-named author's doctoral thesis \cite{palzer-diss} under the guidance of the first-named author. 

  \section{Preliminaries and notation}

  We work in the setting of Ref.~\cite{ap-sph}. Since this reference already contains extensive introductory sections on various aspects of supergeometry and superanalysis, we will restrict ourselves here to only the briefest of comments, thus fixing our notation, and referring to \opcit{} for further details. We do not make any claims to originality in these parts; careful references to the original literature are given in \opcit{}.

  Consider the category of $\cplxs$-superspaces: Its objects are pairs $X=(X_0,\sh O_X)$ comprised of a topological space $X_0$ and a sheaf $\sh O_X$ on $X_0$ of supercommutative $\cplxs$-superalgebras with local stalks; its morphisms $\vphi:X\longrightarrow Y$ are pairs $(\vphi_0,\vphi^\sharp)$ consisting of a continuous map $\vphi_0:X_0\longrightarrow Y_0$ and an even unital morphism of $\cplxs$-superalgebra sheaves $\vphi^\sharp:\sh O_Y\longrightarrow(\vphi_0)_*\sh O_X$ where $(\vphi_0)_*$ denotes the direct image. 

  Given any two $\cplxs$-superspaces $X$ and $Y$, we write $y\in_XY$ to denote the fact that $y:X\longrightarrow Y$ is a morphism of $\cplxs$-superspaces. We call $y$ an \Define{$X$-point} of $Y$. If $f:Y\longrightarrow Z$ is another such morphism, then we write $f(y)$ for $f\circ y$. This enables us to view morphisms as functions on $X$-points, and the Yoneda Lemma from category theory states precise conditions under which such functions are in fact morphisms.

  Let $V$ be a finite-dimensional real super-vector space $V_\ev\oplus V_\odd$, together with a compatible complex structure on the odd part $V_\odd$. The \Define{affine superspace} $\aff(V)$ is defined by 
	\[
	  \aff(V)_0\defi V_\ev,\quad\sh O_{\aff(V)}\defi\sh C^\infty_{V_\ev}\otimes_\reals\textstyle\bigwedge_\cplxs(V_\odd^*).
	\]
	Here, $\sh C^\infty_{V_\ev}$ denotes the sheaf of smooth real-valued functions on $V_\ev$, and $\bigwedge_\cplxs(V_\odd^*)$ is the exterior algebra of the complex vector space $V_\odd^*$. Here and in what follows, we denote the homogeneous parts of a given grading over $\ints/2\ints=\{\ev,\odd\}$ by the subscripts ${}_\ev$ (even) and ${}_\odd$ (odd).

	Given a $\cplxs$-superspace $X$, an \Define{open subspace} is one of the form $X|_U\defi(U,\sh O_X|_U)$ for some open subset $U\subseteq X_0$. A $\cplxs$-superspace $X$ is called a \Define{supermanifold} if $X_0$ is Hausdorff and admits an open cover $(U_i)$ \scth for every index $i$, $X|_{U_i}$ is isomorphic to an open subspace of some affine superspace $\aff(V)$. In particular, the underlying space $X_0$ is naturally a manifold. The category of supermanifolds is defined as the full subcategory of the category of $\cplxs$-superspace whose objects are supermanifolds. One customarily calls supermanifolds as defined above ``\emph{cs} manifolds''. However, we eschew this unfortunate appellation here. 

	The category of supermanifolds admits finite products; thus, there is a notion of group objects in this category, and they will be called \Define{Lie supergroups}. To any Lie supergroup $G$, there is assigned the complex Lie superalgebra $\ger g$ of left-invariant vector fields, and an adjoint action $\Ad$ of $G$ on $\ger g$ \cite[Section 2.4]{ap-sph}. This defines in particular a $G_0$-equivariant Lie superalgebra $(G_0,\ger g)$ \scth $\ger g_\ev$ is the complexification of the Lie algebra of $G_0$ and the differential of $\Ad$ coincides with the restriction of the Lie bracket. Such pairs $(G_0,\ger g)$ of a real Lie group and a complex Lie superalgebra are called \Define{supergroup pairs}, and together with the obvious morphisms, they form a category equivalent to the category of Lie supergroups and their morphisms \cite[Corollary 2.10]{ap-sph}.

	There are natural (categorical) notions of Lie supergroup actions, equivariant, and invariant morphisms. Given a Lie supergroup $G$ and a closed Lie subsupergroup $H$ (\ie $H_0$ is closed in $G_0$ and $\ger h$ is a graded subalgebra of $\ger g$), there is a supermanifold $G/H$ which is universal for $H$-invariant morphisms $G\longrightarrow X$. It admits a natural $G$-action. A detailed exposition of the corresponding theory (in a more general context) is given in \cite[Section 4.4]{ahw-orbits}.

	The theory of integration on a supermanifold $X$ is based on the sheaf $\Abs0{\sh Ber}_X$ of \Define{Berezinian densities} \cite[Definition 2.15]{ap-sph}. It is a finite locally free $\sh O_X$-module of rank $1|0$ or $0|1$ when $X$ has dimension $*|q$ where $q$ is even or odd, respectively. If $x=(u,\xi)$ is a local coordinate system, then there is an associated basis $\Abs0{D(x)}$ of sections. Given a left inverse $r:X\longrightarrow X_0$ to the canonical morphism $X_0\longrightarrow X$ (such an $r$ is called a \Define{retraction}), there is a canonical notion of integrability with respect to $r$ and of the integral $\int_X^r\omega$ of an $r$-integrable section $\omega$ of $\Abs0{\sh Ber}_X$. If $\omega$ has compact support, then it is integrable, and its integral does not depend on $r$. Compare \cite[Section 2.6]{ap-sph} for more details. 

  \section{A recursion formula for spherical superfunctions}\label{sec:recursion}

  In this section, we consider the spherical superfunctions for the super-hyperboloid. We establish a recursive formula in terms of the parameter $\vrho=\tfrac p2-q$, the half sum of positive restricted roots. 

  \subsection{Basic setup}

  We let $(G,K)$ be the pair of Lie supergroups
  \[
    \Parens1{\mathrm{SOUSp}_0(1,1+p\mid 2q),\mathrm{SOSp}(1+p\mid 2q)},\quad p>0,
  \]
  and denote the corresponding pair of Lie superalgebras by $(\ger g,\ger k)$. Here, $G$ and $K$ are determined up to canonical isomorphism by 
  \[
  	G_0\defi\mathrm{SO}_0(1,1+p)\times\mathrm{USp}(2q),\quad K_0\defi\mathrm{SO}(1+p)\times\mathrm{USp}(2q)
  \]
  and 
  \[
  	\ger g\defi\ger{osp}(2+p|2q,\cplxs),\quad\ger k\defi\ger{osp}(1+p|2q,\cplxs),
  \]
  together with the conjugation action of $G_0$ and $K_0$, respectively, on $\ger g$ and $\ger k$. Compare \cite[Section 4.2.2]{ap-sph} for more details. (Note that $\ger{osp}(2+p|2q,\cplxs)$ has no real form whose even part is $\ger{so}(1,1+p)\times\ger{usp}(2q)$, so that we are obliged to abandon the setting of real Lie supergroups and work instead in the present setting.)

  If we define an involution $\theta$ by
  \[
    \theta(x)=\sigma x\sigma,\quad \sigma\defi
      \left(\begin{array}{c:c|c}
          -\mathds{1}_1 &0      &0
        \\\hdashline
          0   &\mathds{1}_{1+p} &0
        \\\hline
          0   &0      &\mathds{1}_q
      \end{array}\right),
  \]
  then $(G,K,\theta)$ is a symmetric supertriple and $(\ger g,\ger k)$ a symmetric superpair in the sense of \cite[Definition 3.1]{ap-sph}.

  The $+1$ eigenspaces of $\theta$ on $\ger g$ is $\ger k$, and its $-1$ eigenspace is denoted by $\ger p$. We let $\ger a\subseteq\ger p_\ev$ be the subspace generated by the matrix
  \begin{equation}\label{eq:h0-def}
    h_0\defi\left(\begin{array}{c:cc|c}
        0 & 1 & 0 & 0
      \\\hdashline
        1 & 0 & 0 & 0
      \\  0 & 0 & 0 & 0
      \\\hline
        0 & 0 & 0 & 0
    \end{array}\right).
  \end{equation}
  We define $\alpha\in\ger a^*$ by $\alpha(h_0)\defi 1$ and frequently identify $\lambda\in\ger a^*$ with the value $\lambda(h_0)$.

  According to \cite[3.1, 4.2.2]{ap-sph}, $\ger g$ decomposes under the action of $\ger a$ as 
  \[
    \ger g=\ger g^{-\alpha}\oplus\ger m\oplus\ger a\oplus\ger g^{\alpha}
  \]
  where 
  \[
    \ger g^{\pm\alpha}\defi\Set1{x\in\ger g}{[h_0,x]=\pm x},\quad\ger m\defi\ger z_\ger k(h_0)=\Set1{x\in\ger k}{[h_0,x]=0}.
  \]
  We also denote $\ger n=\ger g^\alpha$ and $\bar{\ger n}\defi\ger g^{-\alpha}$. Then $\ger k$ and $\ger a$ are $\theta$-invariant, whereas $\bar{\ger n}=\theta(\ger n)$.

  If we let $A$ be the closed subgroup of $G_0$ generated by $\ger a_\reals\defi \reals h_0$ and $N$ the closed connected Lie subsupergroup of $G$ generated by $\ger n$, then the \emph{Iwasawa decomposition} exists, \ie the multiplication morphism 
  \[
    K\times A\times N\longrightarrow G
  \]
  is an isomorphism of supermanifolds \cite[Proposition 3.6]{ap-sph}. We may therefore define morphisms 
  \begin{equation}\label{eq:kan}
    k:G\longrightarrow K,\quad H:G\longrightarrow\aff(\ger a_\reals),\quad n:G\longrightarrow N
  \end{equation}
  by requiring that 
  \begin{equation}\label{eq:iwasawa}
    g=k(g)e^{H(g)}n(g)    
  \end{equation}
  \fa supermanifolds $T$ and all $g\in_T G$.

  We let $\vrho\defi\frac12m_\alpha\cdot\alpha$ where $m_\alpha=\dim\ger n_\ev-\dim\ger n_\odd=p-2q$, the \emph{half sum of positive restricted roots}. It will turn out that the symmetric superfunctions $\phi_\lambda=\phi_\lambda^\vrho$ defined below depend only upon $\lambda$ and $\vrho$.

  Let $M\subseteq K$ be the unique closed Lie subsupergroup with underlying Lie group 
  \[
    M_0\defi Z_{K_0}(h_0)=\Set1{k\in K_0}{kh_0k^{-1}=h_0},
  \]
  and whose Lie superalgebra is $\ger m$. The homogeneous supermanifold $K/M$ carries \cite{ah-ber} an up to multiples unique non-zero Berezinian density $\Abs0{D\dot k}$ which is left $K$-invariant \cite[Definition 3.10]{ap-sph}. We take it to be normalized by
  \begin{equation}\label{eq:km-int-norm}
    \int_{K/M}\Abs0{D\dot k}=\frac{2^{-\vrho}}{\Gamma\Parens1{\vrho+\tfrac12}},
  \end{equation}
  where, as before, we identify $\vrho$ with $\vrho(h_0)$. Such a normalization is possible, since $K/M$ is a supersphere of superdimension $2\vrho=p-2q$, so that \cite[Lemma 4.7]{groeger} applies. The normalization fixes $\Abs0{D\dot k}$ if $\vrho\notin-\tfrac12-\nats$; it does not fix $\Abs0{D\dot k}$ in the other cases.

  We call the supermanifold $G/K$ the \Define{super-hyperboloid}. We now define $K$-invariant superfunctions $\phi_\lambda=\phi_\lambda^\vrho\in\Gamma(\sh O_{G/K})=\Gamma(\sh O_G)^K$, $\lambda\in\ger a^*$, on this supermanifold by 
  \begin{equation}\label{eq:sph-def}
    \phi_\lambda(g)=\phi_\lambda^\vrho(g)\defi\int_{K/M}\Abs0{D\dot k}\,e^{(\lambda-\vrho)(H(gk))}
  \end{equation}
  for all supermanifolds $T$ and all $g\in_TG$. These are the \Define{spherical superfunctions}. Notice that the value of $\phi_\lambda=\phi_\lambda^\vrho$ at unity $1_G$ is the number in Equation \eqref{eq:km-int-norm}.

  The spherical superfunctions $\phi_\lambda=\phi_\lambda^\vrho$ are characterised as follows: They are $K$-biinvariant, so determined uniquely by their restriction to $A$. Here, note that $G$ admits a $KAK$-decomposition in view of \cite[Proposition 3.6]{ap-sph}. The restriction $\phi_\lambda(t)\defi\phi_\lambda(e^{th_0})$ of $\phi_\lambda$ to $A$ satisfies
  \begin{equation}\label{eq:eigenfunction}
      \Lambda_\vrho(\phi_\lambda)=(\lambda^2-\vrho^2)\phi_\lambda,\quad\Lambda_\vrho\defi \partial_t^2+m_\alpha\coth(t)\,\partial_t=\partial_t^2+2\vrho\coth(t)\,\partial_t,
  \end{equation}
  is Weyl-invariant, \ie invariant under $t\longmapsto-t$, and if $\lambda\neq\pm\vrho$, it is up to multiples unique with this property \cite[Proposition 5.1, Theorem 5.5 and proof]{ap-sph}.

  From Equation \eqref{eq:eigenfunction}, it is clear that for $\lambda\neq\pm\vrho$, the function $\phi_\lambda^\vrho(t)$ depends only on $t$, $\lambda$, and $\vrho=\frac p2-q$, but not on $p$ and $q$ separately. Notice that we have not yet fixed the normalization of $\phi_\lambda^\vrho$ for $\vrho\in-\nats-\tfrac12$ and $\lambda\neq\pm\vrho$. We will do this presently, by induction on $\vrho$.

  \subsection{Recursion for the spherical superfunctions}\label{subs:sph-rec}

  Define 
  \[
    \slashed{\partial}\defi\frac1{\sinh t}\partial_t
  \]
  as a differential operator on $\reals\setminus\{0\}$. Then 
  \[
    \sinh^2(t)\,\slashed{\partial}^2=\sinh(t)\partial_t\circ\frac1{\sinh(t)}\partial_t=-\coth(t)\,\partial_t+\partial_t^2,
  \]
  so that 
  \begin{equation}\label{eq:radlap-geod}
    \Lambda_\vrho=\sinh^2(t)\,\slashed{\partial}^2+(2\vrho+1)\cosh(t)\,\slashed{\partial}
  \end{equation}
  where $t$ denotes the identity of $\reals$. Notice that $\slashed{\partial}$ commutes with the Weyl group action $t\longmapsto-t$, so that it maps Weyl-invariant functions to such. In fact, the value $\slashed\partial f(0)$ is well-defined whenever $f$ is Weyl-invariant. 

  \begin{Lem}
    We have 
    \begin{equation}\label{eq:radlap-comm}
      \slashed\partial\Lambda_\vrho-\Lambda_{\vrho+1}\slashed\partial=(2\vrho+1)\slashed\partial.
    \end{equation}
  \end{Lem}

  
  \begin{Prop}[spherical-shift]
    The function $\slashed\partial\phi_\lambda^\vrho$ is a Weyl-invariant eigenfunction of $\Lambda_{\vrho+1}$ to the eigenvalue $\lambda^2-(\vrho+1)^2$. If $\vrho\notin-\nats-\tfrac12$ or $\lambda=\pm\vrho$, then the normalization of $\phi_\lambda^\vrho$ and $\phi_\lambda^{\vrho+1}$ has been fixed \scth 
    \begin{gather}
      \phi_\lambda^\vrho(0)=2^{-\vrho}\Gamma\Parens1{\vrho+\tfrac12}^{-1},\quad \phi_\lambda^{\vrho+1}(0)=2^{-\vrho-1}\Gamma\Parens1{\vrho+\tfrac32}^{-1}\label{eq:spherical-zero},\\
      \slashed\partial\phi_\lambda^\vrho=(\lambda^2-\vrho^2)\phi_\lambda^{\vrho+1}.\label{eq:spherical-shift}
    \end{gather}
    In case $\vrho\in-\nats-\tfrac12$ and $\lambda\neq\pm\vrho$, we may normalize $\phi_\lambda^\vrho$ arbitrarily to achieve Equation \eqref{eq:spherical-zero}. Any choice of normalization for $\phi_\lambda^\vrho$ determines a unique normalization of $\smash{\phi_\lambda^{\vrho+1}}$ ensuring the validity of Equation \eqref{eq:spherical-shift}, and \emph{vice versa}.
  \end{Prop}

  \begin{proof}
    In case $\lambda=\vrho$, $\phi_\lambda^\vrho$ is by definition constant. Since $\phi_\lambda^\vrho(g)=\phi_{-\lambda}^\vrho(g^{-1})$ by \cite[Corollary 3.20]{ap-sph}, the same holds true if $\lambda=-\vrho$. 

    Thus, we may assume that $\lambda\neq\pm\vrho$. By Equation \eqref{eq:radlap-comm}, we find 
    \[
      \Lambda_{\vrho+1}\slashed\partial\phi_\lambda^\vrho=\slashed\partial\Lambda_\vrho\phi_\lambda^\vrho-(2\vrho+1)\slashed\partial\phi_\lambda^\vrho=(\lambda^2-(\vrho+1)^2)\slashed\partial\phi_\lambda^\vrho.
    \]
    This proves the first assertion. 

    It follows that $\slashed\partial\phi_\lambda^\vrho=c\phi_\lambda^{\vrho+1}$ for some constant $c\equiv c(\lambda,\vrho)\in\cplxs$. In particular, 
    \[
      (\lambda^2-\vrho^2)\phi_\lambda^\vrho(0)=(\Lambda_\vrho\phi_\lambda^\vrho)(0)=(2\vrho+1)(\slashed\partial\phi_\lambda^\vrho)(0)=c\,(2\vrho+1)\phi_\lambda^{\vrho+1}(0).
    \]
    If $\smash{\vrho\notin-\tfrac12-\nats}$, then we conclude $c=\lambda^2-\vrho^2$.
    In the remaining cases, if we have chosen a normalization for $\smash{\phi_\lambda^{\vrho+1}}$ or $\smash{\phi_\lambda^\vrho}$, we may still arrange the normalization of $\phi_\lambda^\vrho$ or $\smash{\phi_\lambda^{\vrho+1}}$, respectively, in a unique fashion to achieve $(2\vrho+1)c=\lambda^2-\vrho^2$. 
  \end{proof}

  We have now fixed the normalization of $\phi_\lambda^\vrho$.

  \begin{Cor}
    For $\lambda\neq\pm\vrho$, we have
    \begin{equation}\label{eq:spherical-rec}
      \phi_\lambda^\vrho(t)=(2\vrho+1)\cosh(t)\phi_\lambda^{\vrho+1}(t)+(\lambda^2-(\vrho+1)^2)\sinh^2(t)\phi_\lambda^{\vrho+2}(t).      
    \end{equation}
  \end{Cor}

  \begin{proof}
    This is immediate from Equations \eqref{eq:radlap-geod} and \eqref{eq:spherical-shift}.
  \end{proof}

  \subsection{Recursion for the Harish-Chandra c-function}\label{subs:cfn-rec}

  \begin{Def}
    The \Define{Harish-Chandra $c$-function} is by definition
    \begin{equation}\label{eq:cfn-def}
      \mathbf c(\lambda)=\mathbf c_\vrho(\lambda)\defi\lim\nolimits_{t\to\infty}e^{-(\lambda-\vrho)t}\phi_\lambda^\vrho(t),
    \end{equation}
    whenever this limit exists. 
  \end{Def}

  As shown in \cite[Theorem 4.14]{ap-sph}, the limit $\mathbf c_\vrho(\lambda)$ exists for $\Re\lambda>0$, and the function $\mathbf c(\lambda)$ admits a meromorphic extension to $\ger a^*$ which can be written out explicitly in a super-generalization of the Gindikin--Karpelevi\v c formula. 

  In this subsection, we will not take recourse to this result, but rather only to the classical result that this holds true if $\vrho\sge0$ \cite[Chapter IV, Theorem 6.4]{helgason-gga}. From this, we will give a pedestrian derivation of the general statement, based on the recursion formula established in the previous subsection.

  \begin{Prop}[cfn-recursion]
    If $\Re\lambda>0$ and $\lambda\notin\pm\vrho+\ints$, then $\mathbf c_\vrho(\lambda)$ exists and 
    \begin{equation}\label{eq:cfn-recursion}      
      \mathbf c_{\vrho+1}(\lambda)=\frac2{\lambda+\vrho}\mathbf c_\vrho(\lambda).
    \end{equation}
  \end{Prop}

  \begin{proof}
    By \cite[Chapter IV, Theorem 6.4]{helgason-gga}, the statement holds for $\vrho\sge0$. Assume it holds for $\vrho'\in\vrho+1+\nats$. Then Equation \eqref{eq:spherical-rec} implies that 
    \[
      \begin{split}
        e^{-(\lambda-\vrho)t}\phi_\lambda^\vrho(t)&=\tfrac{2\vrho+1}2(1+e^{-2t})e^{-(\lambda-\vrho-1)t}\phi_\lambda^{\vrho+1}(t)+\tfrac{\lambda^2-(\vrho+1)^2}4(1-e^{-2t})^2\phi_\lambda^{\vrho+2}(t)\\
        &\longrightarrow\tfrac{2\vrho+1}2\mathbf c_{\vrho+1}(\lambda)+\tfrac{\lambda^2-(\vrho+1)^2}4\mathbf c_{\vrho+2}(\lambda)\quad(t\to\infty)
      \end{split}
    \]
    Thus, $\mathbf c_\vrho(\lambda)$ exists and 
    \begin{align*}
      \mathbf c_\vrho(\lambda)&=\tfrac{2\vrho+1}2\mathbf c_{\vrho+1}(\lambda)+\tfrac{\lambda^2-(\vrho+1)^2}4\mathbf c_{\vrho+2}(\lambda).
      \intertext{Applying the inductive assumption to the second summand, this equals}
      &=\tfrac12\Parens1{2\vrho+1+\lambda-\vrho-1}\mathbf c_{\vrho+1}(\lambda)=\frac{\lambda+\vrho}2\mathbf c_{\vrho+1}(\lambda).
    \end{align*}
    This proves the assertion.
  \end{proof}

  The following reproves \cite[Theorem 4.14]{ap-sph}. The case of $\vrho>0$ is a special case of the classical Gindikin--Karpelevi\v{c} formula \cite[Chapter IV, Theorem 6.14]{helgason-gga}.

  \begin{Cor}[gk-fmla]
    If $\Re\lambda>0$ is such that $\lambda\notin\tfrac12\ints$, then, for a suitable choice of normalization of $\phi_\lambda^\vrho$ for $\vrho=\frac12$ and $\lambda\neq\pm\tfrac12$, we have
    \begin{equation}\label{eq:cfn-gkfmla}
      \mathbf c_\vrho(\lambda)=\frac{2^{\vrho-1}\Gamma(\lambda)}{\sqrt\pi\,\Gamma(\lambda+\vrho)}=\frac{2^{-\lambda}\Gamma(\lambda)}{\Gamma\Parens1{\tfrac{\lambda+\vrho}2}\Gamma\Parens1{\frac{\lambda+\vrho+1}2}}.
    \end{equation}
  \end{Cor}

  In the \emph{proof} of the corollary, we need the following lemma.

  \begin{Lem}[sph-cfn-012]
    We have for $\lambda\notin\tfrac12\ints$:
    \begin{equation}\label{eq:sph-012}
      \phi^0_\lambda(t)=\frac1{\sqrt\pi}\cosh(\lambda t),\quad\smash{\phi^{\frac12}_\lambda(t)}=\frac1{\sqrt2}P_{\smash{\lambda-\frac12}}(\cosh(t)),
    \end{equation}
    where $P_s$ denotes the Legendre function. In particular,
    \begin{equation}\label{eq:cfn-012}
      \mathbf c_0(\lambda)=\frac1{2\sqrt\pi},\quad\mathbf c_{\smash{\tfrac12}}(\lambda)=\frac{\Gamma(\lambda)}{\sqrt{2\pi}\,\Gamma\Parens1{\lambda+\tfrac12}}.
    \end{equation}
  \end{Lem}

  \begin{proof}
    The above functions are readily verified to fulfil the eigenfunction relation \eqref{eq:eigenfunction} and to be Weyl-invariant. Moreover, they have the correct values at $0$ since $P_s(1)=1$. We immediately obtain the value of $\mathbf c_0(\lambda)$. For $\mathbf c_{\smash{\tfrac12}}(\lambda)$, recall the formula:
    \[
      P_s(x)=\frac1\pi\int_0^\pi\Parens1{x+\sqrt{x^2-1}\cos\theta}^s\,d\theta,\quad\re s>0,
    \]
    see \cite[Chapter 3.7, (14)]{emot}. This readily implies the claim.
  \end{proof}

  \begin{Rem}
    The expression for $\phi_\lambda^{\frac12}(t)$ in Equation \eqref{eq:sph-012} is to be found in \cite[Chapter IV, Proposition 2.9]{helgason-gga}, for the case of $p=1$, $q=0$. 
  \end{Rem}

  \begin{proof}[\prfof{Cor}{gk-fmla}]
    Recall the classical duplication formula
    \[
      \Gamma(2z)=\frac{2^{2z-1}}{\sqrt\pi}\Gamma(z)\Gamma\Parens1{z+\tfrac12}.
    \]
    \thmref{Lem}{sph-cfn-012} implies the claim for $\vrho=0,\frac12$. The general case follows by induction. 
  \end{proof}

  The following is immediate from \thmref{Prop}{spherical-shift} and \thmref{Prop}{cfn-recursion}.

  \begin{Cor}[sph-der]
    We have 
    \begin{equation}\label{eq:sph-der}
      \slashed\partial\Bracks3{\frac{\phi_\lambda^\vrho}{\mathbf c_\vrho(\lambda)\mathbf c_\vrho(-\lambda)}}(t)=-4\frac{\phi^{\vrho+1}_\lambda(t)}{\mathbf c_{\vrho+1}(\lambda)\mathbf c_{\vrho+1}(-\lambda)}.
    \end{equation}
  \end{Cor}

  \section{The wave packet transform}\label{sec:wavepacket}

  Let $\vphi\in\Gamma(\sh O_{\aff(i\ger a_\reals^*)\times K/M})$. We define the \Define{wave packet transform} by
  \begin{equation}\label{eq:wavepacket-def}
    (\sh J(\vphi))(g)\defi\int_{i\ger a_\reals^*}\frac{d\lambda}{\Abs0{\mathbf c(\lambda)}^2}\int_{K/M}\Abs0{D\dot k}\,\vphi(\lambda,k)e^{-(\lambda+\vrho)(H(g^{-1}k))}
  \end{equation}
  for all supermanifolds $T$ and all $g\in_TG$, provided that the iterated integrals exist. (Note that we use the Iwasawa $KAN$ $A$-projection $H(-)$ instead of the $NAK$ $A$-projection $A(-)$ that Helgason applies. Passage from one to other introduces a shift by $2\vrho$, see \cite[p.~198, Equation (3)]{helgason-sym}.)

  Here and in what follows, we adhere to the convention that $d\lambda$ is the (positive) Hausdorff measure on $i\ger a^*$, \ie $\int_{i\ger a^*}d\lambda$ is $-i$ times the contour integral $\smash{\int_{-i\infty}^{i\infty}dz}$. 

  The result of $\sh J$ is $K$-invariant, and may thus be considered as a function of $\dot g=\pi(g)\in_TG/K$ (where $\pi:G\longrightarrow G/K$ is the canonical projection).

  Notice that if $\vphi(\lambda,k)\equiv\vphi(\lambda)$ is $K$-invariant, then 
  \begin{equation}\label{eq:wavepacket-kinv}
    (\sh J(\vphi))(g)=\int_{i\ger a_\reals^*}\frac{d\lambda}{\Abs0{\mathbf c(\lambda)}^2}\,\phi_\lambda(g)\vphi(\lambda).
  \end{equation}
  Here, we have applied \cite[Corollary 3.20]{ap-sph}.

  \begin{Def}[pw-space][Paley--Wiener space]
    For $k\in\ints$ and $R>0$, define
    \[
      \Norm0{\vphi}_{k,R}\defi\sup\nolimits_{\lambda\in\ger a^*}(1+\Abs0\lambda)^ke^{-R\Abs0{\Re\lambda}}\Abs0{\vphi(\lambda)}
    \]
    for any $\vphi:\ger a^*\longrightarrow\cplxs$. We let $\PW[_R]0{\ger a^*}$ be the space of all entire functions $\vphi:\ger a^*\longrightarrow\cplxs$ \scth $\vphi(\lambda)=\vphi(-\lambda)$ \fa $\lambda\in\ger a^*$ and \fa $k\in\nats$, we have
    \[
      \Norm0{\vphi}_{k,R}<\infty.
    \]
    The inductive limit $\PW0{\ger a^*}\defi\varinjlim_{R>0}\PW[_R]0{\ger a^*}$ of locally convex spaces is called the \Define{Paley--Wiener space} on $\ger a^*$. 
  \end{Def}

  \subsection{Estimates for the derivatives of spherical superfunctions}

  It is clear that $\sh J(\vphi)$ exists for $\vphi\in\PW0{\ger a^*}$. We will improve on this for $\vrho\sle0$ by providing more precise estimates on the growth of $\phi_\lambda$. 

  In this subsection, we assume that $\vrho\sle 0$.

  \begin{Prop}[sphfn-growth]
    For any $\delta>0$ and $\lambda_0\sge0$, there is a constant $C>0$ \scth
    \begin{equation}\label{eq:sphfn-der-est}
      \Abs1{\partial_t^k\phi_\lambda(e^{th_0})}\sle C(1+\Abs0\lambda)^{-\vrho+k}\cosh(\Re\lambda)\cosh^{-\vrho}(t)
    \end{equation}
    \fa $t\in\reals$, $k\sle-\vrho$, and all $\lambda\in\ger a^*\equiv\cplxs$ \scth $\Re\lambda\in[-\lambda_0,\lambda_0]$ and $\Abs0{\lambda-n}>\delta$ \fa $n\in\ints\setminus\{0\}$.
  \end{Prop}

  The \emph{proof} of the proposition is preceded by three lemmas. To state the first, recall the \Define{Harish-Chandra series} 
  \begin{equation}\label{eq:hcseries-def}
    \Phi_\lambda(e^{th_0})\defi e^{(\lambda-\vrho)t}\sum_{\ell=0}^\infty\gamma_\ell(\lambda)e^{-2\ell t}    
  \end{equation}
  where 
  \begin{equation}\label{eq:hcseries-coeff-def}
    \gamma_\ell(\lambda)\defi -\mathbf c(\lambda)\mathbf c(-\lambda)(-1)^\ell\binom{-\vrho}\ell\frac\lambda{(\ell-\lambda)\mathbf c(\ell-\lambda)},
  \end{equation}
  which by \cite[Proposition 5.2]{ap-sph} converges absolutely and uniformly for $t\sge\eps>0$, provided that $\lambda\notin\tfrac12\ints$. Moreover, by \cite[Theorem 5.5]{ap-sph}, we have
  \begin{equation}\label{eq:hcseries-expansion}
    \phi_\lambda(e^{th_0})=(\Phi(\lambda)+\Phi(-\lambda))(e^{th_0})\quad\forall t>0.
  \end{equation}

  \begin{Lem}[hcseries-est]
    Let $\delta>0$ and $\lambda_0\sge0$. There is a constant $C_0>0$ \scth
    \begin{equation}\label{eq:hcseries-estimate}
      \Abs0{\Phi_\lambda(e^{th_0})}\sle C_0\Abs0{\mathbf c(\lambda)}e^{(\Re\lambda-\vrho)t}
    \end{equation}
    \fa $t\sge0$ and all $\lambda\notin\tfrac12\ints$ \scth $\Re\lambda\sle\lambda_0$ and $\Abs0{\lambda-n}\sge\delta$ \fa $n\in\nats\setminus0$.
  \end{Lem}

  \begin{proof}
    We claim that \fa $\ell\in\nats$, we have 
    \begin{equation}\label{eq:hcseries-coeff-estimate}
      \Abs0{\gamma_\ell(\lambda)}\sle c\,\Abs0{\mathbf c(\lambda)}\Abs3{\binom{-\vrho}\ell},
    \end{equation}
    for some positive constant $c$ independent of $\ell$. 

    We prove this claim by induction on $\ell$. To start the induction, assume that $\ell\sge\lambda_0+1-\vrho$ and choose $c$ so that the estimate holds for smaller values of $\ell$. Note that the possibility of doing so depends on the assumption that $\lambda$ is at a positive distance from the set $\nats\setminus0=\{1,2,3,\dotsc\}$ of all non-negative integers. 

    Under the assumption on $\ell$, we have
    \[
      \Abs0{\ell+\vrho-\Re\lambda}=\ell+\vrho-\Re\lambda\sle\ell+1-\Re\lambda,
    \]
    so $\Abs0{\ell+\vrho-\lambda}\sle\Abs0{\ell+1-\lambda}$. This implies
    \[
      \Abs0{\gamma_\ell(\lambda)}=\Abs0{\gamma_{\ell-1}(\lambda)}\,\frac{(\ell-1+\vrho)\Abs0{\ell+\vrho-\lambda}}{\ell\Abs0{\ell+1-\lambda}},
    \]
    which gives Equation \eqref{eq:hcseries-coeff-estimate} by induction.

    Inserting this estimate into Equation \eqref{eq:hcseries-def} yields
    \[
      \Abs0{\Phi_\lambda(e^{th_0})}\sle c\,\Abs0{\mathbf c(\lambda)}e^{(\Re\lambda-\vrho)t}\sum_{\ell=0}^\infty\Abs3{\binom{-\vrho}\ell}e^{-2\ell t}.
    \]
    For $\ell>-\vrho$, we have $\Abs1{\binom{-\vrho}\ell}=\pm(-1)^\ell\binom{-\vrho}\ell$, so 
    \begin{equation}\label{eq:binomabs-trick}
      \sum_{\ell=0}^\infty\Abs3{\binom{-\vrho}\ell}e^{-2\ell t}\sle c'+\Abs3{\sum_{\ell=0}^\infty(-1)^\ell\binom{-\vrho}\ell e^{-2\ell t}}=c'+(1-e^{-2t})^{-\vrho}\sle1+c'
    \end{equation}
    for some positive constant $c'$. This proves the lemma. 
  \end{proof}

  \begin{Lem}[hcseries-der-est]
    Let $\delta>0$ and $\lambda_0\sge0$. For any integer $k\sle-\vrho$, there is a constant $C_k>0$ \scth
    \begin{equation}\label{eq:hcseries-der-estimate}
      \Abs0{\partial_t^k\Phi_\lambda(e^{th_0})}\sle C_k\Abs0{\mathbf c(\lambda)}e^{(\Re\lambda-\vrho)t}
    \end{equation}
    \fa $t\sge0$ and all $\lambda\notin\tfrac12\ints$ \scth $\Re\lambda\sle\lambda_0$ and $\Abs0{\lambda-n}\sge\delta$ \fa $n\in\nats\setminus0$.    
  \end{Lem}

  \begin{proof}
    Notice that 
    \[
      \begin{split}
        \Abs0{\partial_t^ke^{(\lambda-\vrho-2\ell)t}}&\sle\Abs0{\lambda-\vrho-2\ell}^ke^{(\Re\lambda-\vrho-2\ell)t}\\
        &\sle(1+\Abs0\lambda+2\ell-\vrho-1)^ke^{(\Re\lambda-\vrho-2\ell)t}\\
        &\sle c_k\,(1+\Abs0\lambda)^k\ell^ke^{(\Re\lambda-\vrho-2\ell)t},
      \end{split}
    \]
    so that we need to estimate $\sum_\ell\Abs0{\gamma_\ell(\lambda)}\ell^ke^{-2\ell t}$. This can be done along the lines of the proof \thmref{Lem}{hcseries-est}, modifying the estimate in Equation \eqref{eq:binomabs-trick} by inserting a suitable upper bound of $\partial_t^k(1-e^{-2t})^{-\vrho}$.
  \end{proof}

  \begin{Rem}
    Lemmas \ref{Lem:hcseries-est} and \ref{Lem:hcseries-der-est} remain true for $\vrho>0$ and $k>-\vrho$, provided we assume $t>\eps>0$ and allow for $\eps$-dependent constants in Equations \eqref{eq:hcseries-estimate} and \eqref{eq:hcseries-der-estimate}. 

    In case $m_\alpha\sle0$ is even (\ie $-\vrho$ is a non-negative integer), the constants in Equations \eqref{eq:hcseries-estimate} and \eqref{eq:hcseries-der-estimate} become independent of $\delta$ and $\lambda_0$, as the Harish-Chandra series terminates in this case \cite[Corollary 5.4]{ap-sph}.
  \end{Rem}

  \begin{Lem}[cfn-est]
    Let $\delta>0$ and $\lambda_0\sge0$. There is a constant $C_\delta>0$ \scth 
    \begin{equation}\label{eq:cfn-est1}
      \Abs0{\mathbf c(\lambda)}^{-1}\sle C_\delta(1+\Abs0\lambda)^\vrho
    \end{equation}
    \fa $\lambda$ with $\Re\lambda\sge-\lambda_0$ and $\Abs0{\lambda+\vrho+n}\sge\delta$ \fa $n\in\nats$. For $-\vrho$ integral, the estimate holds if we assume only that $\Abs0{\lambda-n}\sge\delta$ for $n=1,\dotsc,-\vrho$.

    We may choose $C_\delta$ \scth in addition
    \begin{equation}\label{eq:cfn-est2}
      \Abs0{\mathbf c(\lambda)}\sle C_\delta(1+\Abs0\lambda)^{-\vrho}
    \end{equation}
    \fa $\lambda$ with $\Re\lambda\sge-\lambda_0$ and $\Abs0{\lambda+\vrho+n}\sge\delta$ \fa $n\in\nats$.
  \end{Lem}

  \begin{proof}
    From Equation \eqref{eq:cfn-gkfmla}, it follows that $\mathbf c(\lambda)$ has but finitely many zeros with $\Re\lambda\sge-\lambda_0$, and they are all of the shape $\lambda\equiv-\vrho-n$ for some non-negative integer $n$. Moreover, we have 
    \[
      \lim\nolimits_{\Abs0\lambda\to\infty}\Parens3{\lambda^\vrho\frac{\Gamma(\lambda)}{\Gamma(\lambda+\vrho)}}=1,
    \]
    provided that $\Abs0{\arg\lambda}<\pi$, so the assertion follows. 
  \end{proof}

  \begin{proof}[\prfof{Prop}{sphfn-growth}]
    For $\lambda\notin\tfrac12\ints$, $\Re\lambda>0$ and $t>0$, the estimate in Equation \eqref{eq:sphfn-der-est} follows from Lemmas \ref{Lem:hcseries-est} and \ref{Lem:hcseries-der-est} upon using Equation \eqref{eq:hcseries-expansion}, together with 
    \[
      \frac{e^t}2\sle\cosh(t)\sle e^t,\quad\Abs0{e^{\lambda t}}\sle2\cosh(\Re\lambda t).
    \]
    The general case follows from the symmetry and continuity of $\phi_\lambda(t)$ as a function of $t$ and $\lambda$.
  \end{proof}

  \begin{Cor}[sphfn-wt-growth]
    Let $\delta>0$, $\lambda_0\sge0$, and $k\sle -\vrho$ a non-negative integer. Then there is a constant $C>0$ \scth 
    \begin{equation}\label{eq:sphfn-wt-der-est}
          \Abs3{\frac{\partial_t^k\phi_\lambda(e^{th_0})}{\mathbf c(\lambda)\mathbf c(-\lambda)}}\sle C(1+\Abs0\lambda)^{\vrho+k}\cosh(\Re\lambda t)\cosh^{-\vrho}(t)
    \end{equation}
    \fa $t\in\reals$, $\lambda$ \scth $\Re\lambda\in[-\lambda_0,\lambda_0]$ and $\Abs0{\lambda+\vrho+n}\sge\delta$ for all $n\in\nats$, $n\neq-\vrho$.
  \end{Cor}

  \begin{proof}
    Similar to the proof of \thmref{Prop}{sphfn-growth}, it will be sufficient to prove the statement for $t>0$, $\Re\lambda>0$, and $\lambda\notin\tfrac12\ints$. Under these assumptions, it will follow from the existence of some positive constant $C'$ independent of $t$ \scth 
    \begin{equation}\label{eq:hcseries-wt-der-est}
          \Abs3{\frac{\partial_t^k\Phi_\lambda(e^{th_0})}{\mathbf c(\lambda)\mathbf c(-\lambda)}}\sle C'(1+\Abs0\lambda)^{\vrho+k}e^{(\Re\lambda-\vrho)t}
    \end{equation}
    \fa $\lambda$ \scth $\Re\lambda\in[-\lambda_0,\lambda_0]$ and $\Abs0{\lambda+\vrho+n}\sge\delta$ for all $n\in\nats$, $n\neq-\vrho$. 

    In turn, the estimate in Equation \eqref{eq:hcseries-wt-der-est} follows from \thmref{Lem}{hcseries-der-est}, its proof, and Equation \eqref{eq:hcseries-coeff-def}. Indeed, Equation \eqref{eq:cfn-est1} implies
    \[
      \Abs0{\mathbf c(\ell-\lambda)}\sle C_\delta(1+\Abs0{\ell-\lambda})^\vrho\sle C_\delta(1+\Abs0\lambda)^\vrho
    \]
    for $\ell>2\lambda_0$, because $\Abs0{\ell-\lambda}^2=\Abs0\lambda^2+\ell(\ell-2\Re\lambda)\sge\Abs0\lambda^2$ in that case.
  \end{proof}

  \subsection{Residues of the wave packet transform}

  In this subsection, we discuss the existence and the asymptotics of the wave packet transform for $\vrho\sle 0$. To that end, let $A^+\defi\Set1{e^{th_0}}{t>0}$ be the (exponential image of the) \Define{positive Weyl chamber}.

  \begin{Lem}[wavepacket-hcseries]
    Let $\vphi:i\ger a_\reals^*\longrightarrow\cplxs$ be a continuous function. \Fa $a\in A$, we have $(\sh J(\vphi))(a^{-1})=(\sh J(\vphi))(a)$ in the sense that the left-hand side of the equation exists if and only so does the right, and in this case, we have the equality.

    Moreover, for $a\in A^+$, we have 
    \begin{equation}\label{eq:wavepacket-hcseries}
      (\sh J(\vphi))(a)=2\int_{i\ger a_\reals}\frac{d\lambda}{\mathbf c(\lambda)\mathbf c(-\lambda)}\Phi_\lambda(a)\vphi(\lambda),
    \end{equation}    
    again in the sense of simultaneous existence of both sides of the equation and equality in case of existence.
  \end{Lem}

  \begin{proof}
    We have $\mathbf c(-\lambda)=\overline{\mathbf c(\lambda)}$ and $\phi_\lambda(a)=\phi_{-\lambda}(a^{-1})=\phi_{-\lambda}(a)$. This implies the first equality. Moreover, Equation \eqref{eq:hcseries-expansion} implies Equation \eqref{eq:wavepacket-hcseries}.
  \end{proof}

  \begin{Prop}
    Assume that $\vrho\sle0$. Let $R>0$ and $n\in\ints$ \scth $\vrho+1\sle n$. If $\Norm0\vphi_{n+1,R}<\infty$, then $\sh J(\vphi)(e^{th_0})$ exists as a function of $t\in\reals$ and has continuous derivatives up to order $-\vrho+n-1$ in $t$.
  \end{Prop}

  \begin{proof}
     The claim follows from \eqref{eq:wavepacket-hcseries} and Equation \eqref{eq:hcseries-wt-der-est}. 
  \end{proof}

  Due to the $K$-biinvariance of $\phi_\lambda$, the following is immediate.

  \begin{Cor}
    If $\vphi\in\PW0{\ger a^*}$, then $\sh J(\vphi)$ exists as a superfunction on $G/K$.
  \end{Cor}

  Moreover, we obtain the following fact which appears counter-intuitive from the classical ungraded theory.

  \begin{Cor}
    If $\vrho<-k-1$ where $k\in\nats$, then $\sh J(1)$ exists on $A$ and is $k$-times continuously differentiable.
  \end{Cor}

  These results can be sharpened somewhat if we take $\vphi$ to be entire.

  \begin{Prop}[wavepacket-res]
    Assume $\vrho\sle0$ and let $R\sge0$. Let $\vphi:\ger a^*\longrightarrow\cplxs$ be entire \scth $\Norm0\vphi_{n,R}<\infty$ \fs integer $n>\vrho$. Then $\sh J(\vphi)(e^{th_0})$ exists as an improper integral for $t>R$, and in this case
    \[
      \sh J(\vphi)(e^{th_0})=4\pi\sum_{\Abs0k<-\vrho}\vphi(\vrho+k)\,\res_{\lambda=\vrho+k}\Bracks3{\frac{\Phi_\lambda(e^{th_0})}{\mathbf c(\lambda)\mathbf c(-\lambda)}}.
    \]
  \end{Prop}

  The \emph{proof} depends on the following lemma.

  \begin{Lem}[wavepacket-residues]
    Let $\vrho\sle0$ and fix $t>0$. The poles of $[\mathbf c(\lambda)\mathbf c(-\lambda)]^{-1}\Phi_\lambda(e^{th_0})$, considered as a meromorphic function of $\lambda$, are located at $\lambda=\vrho+k$ where $k\in\nats$, $k<-\vrho$. The residues have the values
    \begin{align}
      \res_{\lambda=\vrho+k}\Bracks3{\frac{\Phi_\lambda(e^{th_0})}{\mathbf c(\lambda)\mathbf c(-\lambda)}}
      &=\frac{\sqrt\pi\,2^{1-\vrho}(-1)^k(\vrho+k)}{\Gamma(1-\vrho)}\sum_{\ell=0}^\infty\binom{-\vrho}\ell\binom{-\vrho}{k-\ell}e^{(k-2\ell)t}\label{eq:wavepacket-residues1}\\
      &=\frac{\sqrt\pi\,2^{1-\vrho}(\vrho+k)}{\Gamma(1-\vrho)k!}\,\partial^k_{s=0}(1-2s\cosh(t)+s^2)^{-\vrho}\label{eq:wavepacket-residues2}.
    \end{align}
  \end{Lem}

  \begin{proof}
    We expand for small $\Abs0{s}$:
    \[
      \begin{split}
        (1+2s\cosh(t)+s^2)^{-\vrho}
        &=\sum_{n=0}^\infty s^n\sum_{\ell=0}^\infty\binom{-\vrho}\ell\binom{-\vrho}{n-\ell}e^{(n-2\ell)t},
      \end{split}
    \]
    Therefore, Equation \eqref{eq:wavepacket-residues2} will follow from Equation \eqref{eq:wavepacket-residues1}.

    From Equations \eqref{eq:hcseries-coeff-def} and \eqref{eq:cfn-gkfmla}, we have
    \[
      \begin{split}
        \frac{\gamma_\ell(\lambda)}{\mathbf c(\lambda)\mathbf c(-\lambda)}
        &=\sqrt\pi\,2^{1-\vrho}(-1)^{\ell+1}\binom{-\vrho}\ell\frac{\lambda\,\Gamma(\ell+\vrho-\lambda)}{\Gamma(\ell+1-\lambda)}.
      \end{split}
    \]
    The poles of $[\mathbf c(\lambda)\mathbf c(-\lambda)]^{-1}\Phi_\lambda(e^{th_0})$ are located at $\vrho+\nats$. Any poles at zero are eliminated by the factor $\lambda$. The function $\Gamma(\ell+\vrho-\lambda)$ has a pole at $\lambda=\vrho+k$, and its residue is $(-1)^{k+\ell+1}\frac1{(k-\ell)!}$ for $\ell\sle k$ and $0$ otherwise.

    Thus, we obtain non-zero residues only for $\lambda=\vrho+k$, $k<-\vrho$, and they are 
    \[
      \begin{split}
        \res_{\lambda=\vrho+k}\Bracks3{\frac{\gamma_\ell(\lambda)}{\mathbf c(\lambda)\mathbf c(-\lambda)}}
        &=\sqrt\pi\,2^{1-\vrho}\frac{(-1)^k(\vrho+k)}{\Gamma(1-\vrho)}\binom{-\vrho}\ell\binom{-\vrho}{k-\ell}.
      \end{split}
    \]
    This proves Equation \eqref{eq:wavepacket-residues1} by the definition of $\Phi_\lambda$ in Equation \eqref{eq:hcseries-def}.
  \end{proof}

  \begin{Rem}
    Similar computations as in the proof of the previous lemma show that 
    \[
      \phi_\lambda^\vrho(t)={}_2F_1\Parens1{\tfrac{\rho+\lambda}{2},\tfrac{\rho-\lambda}{2};\tfrac12+\vrho;-\sinh^2(t)}.
    \]
  \end{Rem}

  \begin{proof}[\prfof{Prop}{wavepacket-res}]
    We use Equation \eqref{eq:wavepacket-hcseries} to express the wave packet transform. Then this becomes a somewhat standard contour shift argument: Let $\lambda_0>\Abs0\vrho$ and define a contour by $\gamma(\theta)\defi i\lambda_0e^{i\theta}$ for $\theta\in[0,\pi]$. Then by the residue formula
    \begin{equation}\label{eq:res-fmla}
      i\int_{-i\lambda_0}^{i\lambda_0}\frac{d\lambda}{\mathbf c(\lambda)\mathbf c(-\lambda)}\,\Phi_\lambda(e^{th_0})\vphi(\lambda)+\int_\gamma\frac{d\lambda}{\mathbf c(\lambda)\mathbf c(-\lambda)}\,\Phi_\lambda(e^{th_0})\vphi(\lambda)=2\pi iR_\vphi,
    \end{equation}
    where 
    \[
      R_\vphi\defi\sum_{\Abs0k<-\vrho}\vphi(\vrho+k)\,\res_{\lambda=\vrho+k}\Bracks3{\frac{\Phi_\lambda(e^{th_0})}{\mathbf c(\lambda)\mathbf c(-\lambda)}},
    \]
    since by \thmref{Lem}{wavepacket-residues}, the poles of the integrand lie in the interior of the contour $\gamma([-i\lambda_0,i\lambda_0])$ and are of the shape $\lambda\equiv\vrho+k$ for non-negative integers $k<-\vrho$. 

    The numbers $\lambda_0$, $\lambda_0\to\infty$, may be chosen so that we are at a positive distance from the integers. Thus, Equation \eqref{eq:hcseries-wt-der-est} applies, and we find for $\lambda=\gamma(\theta)$:
    \[
      \begin{split}
        \Abs3{\frac1{\mathbf c(\lambda)\mathbf c(-\lambda)}\,\Phi_\lambda(e^{th_0})\vphi(\lambda)}
        &\sle C(1+\Abs0\lambda)^{\vrho-n}e^{\Re\lambda(t-R)-\vrho t}\\
        &= C(1+\lambda_0)^{\vrho-n}e^{-\vrho t}e^{(R-t)\lambda_0\sin\theta}.
      \end{split}
    \]
    Since $\sin\theta$ is symmetric around $\frac\pi2$ and $\sin\theta\sge\frac\theta2$ for $\theta\in[0,\tfrac\pi2]$, we find that 
    \[
      \begin{split}
        \Abs3{\int_\gamma\frac{d\lambda}{\mathbf c(\lambda)\mathbf c(-\lambda)}\,\Phi_\lambda(e^{th_0})\vphi(\lambda)}&\sle
        2C(1+\lambda_0)^{\vrho-n}e^{-\vrho t}\int_0^{\frac\pi2}d\theta\,\lambda_0e^{(R-t)\lambda_0\frac\theta2}\\
        &\sle4C(1+\lambda_0)^{\vrho-n}e^{-\vrho t}\frac1{t-R}.
      \end{split}
    \]
    This quantity vanishes for $\lambda_0\longrightarrow\infty$, so that Equation \eqref{eq:res-fmla} shows the claim. 
  \end{proof}

  Applying \thmref{Prop}{wavepacket-res} for $n=0$ and $R=0$, we obtain the following corollary.

  \begin{Cor}[wavepacket-const]
    Assume that $\vrho<0$. Then $\sh J(1)(e^{th_0})$ exists for $t\neq0$ and equals
    \begin{equation}\label{eq:wavepacket-const}
      \begin{split}
        \sh J(1)(e^{th_0})&=4\pi\sum_{\Abs0k<-\vrho}\res\limits_{\lambda=\vrho+k}\Bracks3{\frac{\Phi_\lambda(e^{th_0})}{\mathbf c(\lambda)\mathbf c(-\lambda)}}\\
        &=-\frac{\pi^{\frac32}\,2^{2-\vrho}}{\Gamma(1-\vrho)\Gamma(-2\vrho)}\,\partial_{s=0}^{-2\vrho-1}\Bracks3{\frac{(1-2s^2\cosh(t)+s^4)^{-\vrho}}{(1-s)^2}}.
      \end{split}
    \end{equation}
    In particular, $\sh J(1)$ admits a smooth extension to $A$ which in what follows will be denoted by the same symbol. 
  \end{Cor}

  \begin{proof}
    It remains only to prove the second equality in Equation \eqref{eq:wavepacket-const}. From \thmref{Prop}{wavepacket-res} and Equation \eqref{eq:wavepacket-residues2}, we find 
    \begin{align*}
      \sh J(1)(e^{th_0})&=\frac{\pi^{\frac32}\,2^{3-\vrho}}{\Gamma(1-\vrho)}\sum_{k<-\vrho}\frac{\vrho+k}{k!}\partial^k_{s=0}(1-2s\cosh(t)+s^2)^{-\vrho}.
      \intertext{Since for smooth $f$, the Taylor series of $f(s^2)$ is $\sum_{k=0}^\infty\frac{s^{2k}}{k!}\partial^k_{s=0}f(s)$, this equals} 
      &=\frac{\pi^{\frac32}\,2^{2-\vrho}}{\Gamma(1-\vrho)}\sum_{k=0}^{-2\vrho-1}\frac{2\vrho+k}{k!}\partial^k_{s=0}(1-2s^2\cosh(t)+s^4)^{-\vrho}.
      \intertext{Moreover, the Taylor coefficient of order~$-2\vrho-k-1$ of $(1-s)^{-2}$ is $-2\vrho-k$, so the Leibniz formula shows that the result is}
      &=-\frac{\pi^{\frac32}\,2^{2-\vrho}}{\Gamma(1-\vrho)\Gamma(-2\vrho)}\,\partial^{-2\vrho-1}_{s=0}\Bracks3{\frac{(1-2s^2\cosh(t)+s^4)^{-\vrho}}{(1-s)^2}},
    \end{align*}
    thereby proving the claim.
  \end{proof}

  \begin{Cor}[wavepacket-const-halfint]
    Assume that $\vrho<0$ is half-integral. Then for $t\neq0$, we have
    \begin{equation}\label{eq:wavepacket-const-halfint}
      \sh J(1)(e^{th_0})=-\frac{\pi\,2^{3-2\vrho}}{\sqrt2\,\Gamma\Parens1{\tfrac12-\vrho}}\,(1-\cosh(t))^{-\vrho-\frac12}.
    \end{equation}
  \end{Cor}

  \begin{proof}
    Since $-2\vrho-1=-2(\vrho+\tfrac12)$ is an even integer, $\partial^{-2\vrho-1}_{s=0}$ vanishes on odd functions. Decomposing in Equation \eqref{eq:wavepacket-const}
    \[
      \frac1{(1-s)^2}=\frac{1+s^2}{(1-s^2)^2}+\frac{2s}{(1-s^2)^2},
    \]
    only the first summand contributes in the derivative. Again using the fact that the Taylor series of $f(s^2)$ equals $\sum_{k=0}^\infty\smash{\frac{s^{2k}}{k!}}\partial^k_{s=0}f(s)$, we find that
    \begin{equation}\label{eq:wavpacket-const-analytic}
      \sh J(1)(e^{th_0})=\frac{-\pi^{\frac32}\,2^{2-\vrho}}{\Gamma(1-\vrho)\Gamma\Parens1{\tfrac12-\vrho}}\,\partial^{-\vrho-\tfrac12}_{s=0}\Bracks3{\frac{(1+s)(1-2s\cosh(t)+s^2)^{-\vrho}}{(1-s)^2}}.      
    \end{equation}
    The $k$th Taylor coefficient of $\Parens1{\Parens1{-\vrho-\tfrac12}!}^{-1}\partial^{-\vrho-\tfrac12}_{s=0}\Bracks1{\frac{(1+s)(1-2sx+s^2)^{-\vrho}}{(1-s)^2}}$ at $x=1$ is 
    \begin{align*}
      (-2)^k\binom{-\vrho}k&\frac{\partial^{-\vrho-\tfrac12}_{s=0}}{\Parens1{-\vrho-\tfrac12}!}\Bracks3{s^k\frac{(1+s)(1-s)^{-2(\vrho+k)}}{(1-s)^2}}\\
      &=(-2)^k\binom{-\vrho}k\frac{\partial^{-\vrho-\tfrac12}_{s=0}}{\Parens1{-\vrho-\tfrac12}!}\Parens1{s^k(1+s)(1-s)^{-2(\vrho+k+1)}}.
      \intertext{This is non-zero only if $k\sle -\vrho-\tfrac12$. If $k\sle-\vrho-\tfrac32$, it equals} 
      &=(-2)^k\binom{-\vrho}k\Bracks3{\frac{\partial^{-\vrho-\tfrac12-k}_{s=0}}{\Parens1{-\vrho-\tfrac12-k}!}+\frac{\partial_{s=0}^{-\vrho-\tfrac32-k}}{\Parens1{-\vrho-\tfrac32-k}!}}(1-s)^{-2(\vrho+k+1)}\\
      &=(-2)^k\binom{-\vrho}k\Bracks3{(-1)^{m}\binom{2m-1}m+(-1)^{m-1}\binom{2m-1}{m-1}}=0,\\
      \intertext{for $m\defi-\vrho-\tfrac12-k$, and for $k=-\vrho-\tfrac12$, it is}
      &=(-2)^{-\vrho-\tfrac12}\binom{-\vrho}{-\vrho-\tfrac12}\frac{1+s}{1-s}\Big|_{s=0}=(-2)^{-\vrho-\frac12}\binom{-\vrho}{-\vrho-\tfrac12}.
    \end{align*}
    Since $\sh J(1)(e^{th_0})$ is analytic in $t$ by inspection of Equation \eqref{eq:wavpacket-const-analytic}, we may expand to obtain
    \[
      \begin{split}
        \sh J(1)(e^{th_0})
        &=-\frac{\pi\,2^{3-2\vrho}}{\sqrt2\,\Gamma\Parens1{\tfrac12-\vrho}}\,(1-\cosh(t))^{-\vrho-\frac12},
      \end{split}
    \]
    as was claimed.
  \end{proof}

  \section{The inversion formula for the Helgason--Fourier transform}\label{sec:inversion}

  In this section, we turn to our main objective, the proof of an inversion formula for the non-Euclidean Fourier transform on the super-hyperbolic space $G/K$.

  \subsection{Invariant Berezin integration in polar coordinates}

  In this subsection, we give explicit expressions for the invariant Berezinian integrals on $G/K$ and $K/M$, generalizing the `polar coordinates' on a Riemannian symmetric space of non-compact type. We begin by presenting explicit models for these supermanifolds. 

  For any supermanifold $T$ and any $g\in_TG$, we decompose
  \[
    g=
    \left(
      \begin{array}{c:c}
        a&b\\
        \hdashline
        c&d
      \end{array}
    \right)
  \]
  where $a$ is a $1\times 1$ matrix, $b$ is a $1\times(1+p|2q)$ matrix, $c$ is a $(1+p|2q)\times1$ matrix, and $d$ is a $(1+p|2q)\times(1+p|2q)$ matrix. 

  We let $D$ be the open subspace of $\aff^{1+p|2q}$ whose underlying open set is the Euclidean unit ball of $\reals^{1+p}$. Then $g$ acts on $x\in_TD$ by 
  \[
    g\cdot x\defi(dx+c)(bx+a)^{-1}.
  \]
  This defines an action of the Lie supergroup $G$ on the supermanifold $D$.

  Taking $o\defi 0\in D_0$ as our base point, we compute that the isotropy supergroup of the action at $o$ to be the closed subsupergroup representing the functor
  \[
    G_o(T)=\Set3{g=\left(
      \begin{array}{c:c}
        a&b\\
        \hdashline
        c&d
      \end{array}
    \right)\in_TG}{b=c=0,a=1}=K(T),
  \]
  see \cite[Theorem 4.20]{ahw-orbits}. Therefore, we have a $G$-equivariant injective immersion $G/K\longrightarrow D$ mapping $K_0$ to $o$ \cite[Theorem 4.24]{ahw-orbits}. Since $\dim G/K=1+p|2q=\dim D$, this immersion is an open embedding. Since the underlying action is transitive, it follows that $G/K$ is $G$-equivariantly isomorphic to $D$. 

  Similarly, $K$ acts linearly on $\aff^{1+p|2q}$, leaving invariant the closed subsupermanifold $B$ of $\aff^{1+p|2q}$ defined by the even equation 
  \[
    \Norm0x^2\defi\sum_{j=1}^p(x^j)^2+2\sum_{j=1}^qx^{p+2j-1}x^{p+2j}=1,
  \]
  \ie the supersphere of dimension $p|2q$. The isotropy subsupergroup at the point $e_1\defi(1,0,\dotsc,0)^t\in B_0$ is precisely $M$, and $\dim K/M=p|2q=\dim B$, so as above, $K/M$ is $K$-equivariantly isomorphic to $B$.

  In the following, we shall identify $G/K$ with $D$ and $K/M$ with $B$. We will use the following \Define{$KA$ decomposition}
  \begin{equation}\label{eq:kak-decomp}
    ka\cdot o=\tanh(t)(k\cdot e_1),\quad a=e^{th_0}\in_TA^+,    
  \end{equation}
  for any supermanifold $T$, any $t\in\Gamma(\sh O_{T,\ev})$, $t_0>0$, and any $k\in_TK$. We use it to define $r\defi\tanh(t)$, where $t$ is the identity of $(0,\infty)$, and 
  \begin{equation}\label{eq:aplus-param}
    a_r\defi e^{th_0}:(0,1)\longrightarrow A^+.
  \end{equation}

  The metric given by polarization of $\Norm0\cdot^2$ on $\aff^{1+p|2q}$ induces a Riemannian metric on $B$. The Riemannian Berezinian density on $B$ (\cf \cite{groeger}) will be denoted by $\Abs0{Db}$. Moreover, $D$ carries a non-zero $G$-invariant Berezinian density $\Abs0{D\dot g}$ which is unique up to multiples \cite{ah-ber}. We presently fix its normalization by relating it to the Riemannian (\ie Lebesgue) Berezinian density $\Abs0{D(x)}$ on $\aff^{1+p|2q}$.

  Notice that by Equation \eqref{eq:km-int-norm} and the normalization of the spherical superfunctions $\phi_\lambda^\vrho$ in Subsections \ref{subs:sph-rec} and \ref{subs:cfn-rec}, we have fixed the normalization of $\Abs0{D\dot k}$. The normalization is easy to determine in case $\vrho\notin-\nats-\tfrac12$. 

  \begin{Lem}[gk-km-ber]
    The density $\Abs0{Db}$ on $B$ is $K$-invariant, so that $\Abs0{D\dot k}=c_{K/M}\Abs0{Db}$ for some non-zero constant $c_{K/M}$. If $\vrho\notin-\nats-\tfrac12$, then this constant is 
    \begin{equation}\label{eq:km-ber}
      c_{K/M}=\frac{(-1)^q}{\sqrt2(2\pi)^{\frac{p+1}2}}.
    \end{equation}
    The density $\Abs0{D(x)}\,(1-\Norm0x^2)^{-1-\vrho}$ on $D$, where $x$ is the identity of $D$, is $G$-invariant, so that we may normalize $\Abs0{D\dot g}$ by 
    \begin{equation}\label{eq:gk-ber}
      \Abs0{D\dot g}\defi c_{K/M}\cdot\frac{\Abs0{D(x)}}{(1-\Norm0x^2)^{1+\vrho}}. 
    \end{equation}
  \end{Lem}

  \begin{proof}
    The $K$-invariance of $\Abs0{Db}$ is easy to check, since the action in question is linear. Then Equation \eqref{eq:km-ber} follows from Equation \eqref{eq:km-int-norm} and \cite[Lemma 4.7]{groeger}. We record the fact that 
    \begin{equation}\label{eq:volb}
      \int_B\Abs0{Db}=(-1)^q\frac{2^{q+1}\pi^{\frac{p+1}2}}{\Gamma\Parens1{\vrho+\frac12}}.
    \end{equation}
    In case $\vrho$ is a negative half-integer, we may take $\Abs0{D\dot k}$ to be any non-zero multiple of $\Abs0{Db}$ to achieve Equation \eqref{eq:km-int-norm}.

    Similarly, $\Abs0{D(x)}\cdot(1-\Norm0x^2)^{-1-\vrho}$ is invariant under the action of $K$. It is thus sufficient to check that it is invariant under the action of $A$. But
    \[
      \ABer3{\frac{\partial(e^{th_0}\cdot x)}{\partial x}}=\Parens1{x^1\sinh(t)+\cosh(t)}^{-2-2\vrho}=\Bracks3{\frac{1-\Norm0{e^{th_0}\cdot x}^2}{1-\Norm0x^2}}^{1+\vrho},
    \]
    where $x=(x^1,\dotsc)$ is the standard coordinate system on $\aff^{1+p|2q}$, as follows from the identity
    \[
      (sx^1+c)^2-(cx^1+s)^2=1-(x^1)^2,\quad c\defi\cosh(t),s\defi\sinh(t).
    \]
    This entails the assertion. 
  \end{proof}

  We can now give the desired explicit expression for the invariant integral on $G/K$. In the following, recall that $2\vrho=p-2q$.

  \begin{Th}[inv-ber-gk]
    Let $f\in\Gamma_c(\sh O_{G/K})=\Gamma_c(\sh O_G)^K$.
    \begin{enumerate}[wide]
      \item\label{item:inv-ber-gk-i} Let $\vrho\sge0$. Then 
      \begin{equation}\label{eq:inv-ber-gk-i}
        \int_{G/K}\Abs0{D\dot g}\,f(g)=\int_0^1dr\,r^{2\vrho}(1-r^2)^{-1-\vrho}\int_B\Abs0{D\dot k}\,f(ka_r).
      \end{equation}
      \item\label{item:inv-ber-gk-ii} Let $\vrho<0$ be integral. Then 
      \begin{equation}\label{eq:inv-ber-gk-ii}
        \int_{G/K}\Abs0{D\dot g}\,f(g)=\frac{(-1)^\vrho\Gamma\Parens1{\vrho+\tfrac12}}{2\sqrt\pi}\!\int_0^1\frac{dr}{\sqrt r}\,\partial_r^{-\vrho}\Bracks3{(1-r)^{-1-\vrho}\!\int_B\Abs0{D\dot k}\,f(ka_{\sqrt r})}.
      \end{equation}
      \item\label{item:inv-ber-gk-iii} Let $\vrho<0$ be half-integral. Then 
      \begin{equation}\label{eq:inv-ber-gk-iii}
        \begin{split}
          \int_{G/K}\Abs0{D\dot g}\,f(g)=\int_0^1&dr\,r^{2\vrho}(1-r^2)^{-1-\vrho}\int_B\Abs0{D\dot k}\,f(ka_r)\\
          +&\frac{\partial_{r=0}^{-2\vrho-1}}{(-2\vrho-1)!}\Bracks3{(1-r^2)^{-1-\vrho}\int_B\Abs0{D\dot k}\,\log\Norm0{ke_1}_\ev f(ka_r)},
        \end{split}
      \end{equation}
      where we define $\Norm0b_\ev\defi\sqrt{(b^1)^2+\dotsm+(b^p)^2}$, with $b$ denoting the identity of $B$.
    \end{enumerate}
  \end{Th}

  \begin{Rem}
    The theorem recovers the integration in polar coordinates for the classical hyperbolic space, see \cite[Chapter I, Theorem 5.8]{helgason-gga}. Moreover, as follows from Equation \eqref{eq:polar-bdterms-rhoneg-halfint} below, it also recovers Zirnbauer's corresponding formula \cite[Theorem 1]{zirn-cmp} for the super-Poincar\'e disc, where $\vrho=-\tfrac12$.
  \end{Rem}

  We will \emph{prove} the theorem case-by-case. First, let 
  \[
    x=(x^1,\dotsc,x^{1+p+2q})=(u,\xi)=(u^1,\dotsc,u^{1+p},\xi^1,\dotsc,\xi^{2q})
  \]
  be the standard coordinate system on $\aff^{1+p|2q}$ and $\gamma$ the standard retraction of $\aff^{1+p|2q}$, determined by $\gamma^\sharp(u_0)=u$. In order to introduce super-polar coordinates, we let $U_0\defi\reals^{1+p}\setminus((-\infty,0]\times\reals^p)$ and denote by $v_0$ the standard polar coordinates on $U_0$ defined as usual by the assumption that the tuple $v_0=(v^1_0,\dotsc,v^{1+p}_0)$ of functions is valued in $V_0\defi(0,\infty)\times(-\pi,\pi)\times\Parens1{-\tfrac\pi2,\tfrac\pi2}^{p-1}$ and 
  \begin{equation}\label{eq:polar}
    u^j_0=
    \begin{cases}
      v^1_0\sin(v^{1+p}_0), &\text{for }j=1+p,\\
      v^1_0\sin(v^j_0)\cos(v^{j+1}_0)\dotsm\cos(v^{1+p}_0), &\text{for }j=2,\dotsc,p,\\
      v^1_0\cos(v^2_0)\cos(v^3_0)\dotsm\cos(v^{1+p}_0), &\text{for }j=1.
    \end{cases}
  \end{equation}
  In particular, we have $R_0\defi\Norm0{u_0}=v^1_0$. 

  We now define a retraction $\gamma'$ of $\aff^{1+p|2q}\setminus\{0\}$ by $R\defi\Norm0x$ and 
  \[
    (\gamma')^\sharp(u_0)=\frac{R}{\gamma^\sharp(R_0)}\,u.
  \]
  We set $\eta\defi R^{-1}\xi$ and follow the convention that $\eta$ is a column while $\eta^*$ is the row with entries $\eta^{*j}\defi\eta^{j-1}$ when $j$ is even and $\eta^{*j}\defi-\eta^{j+1}$ when $j$ is odd, so that
  \[
    R^2=\Norm0u^2+\xi^*\xi=\Norm0u^2+R^2\,\eta^*\eta,
  \]
  and thus 
  \begin{equation}\label{eq:super-radius}
    \gamma^\sharp(R_0)=\Norm0u=R\sqrt{1-\eta^*\eta}.
  \end{equation}
  Defining
  \[
    v=(v^1,\dotsc,v^{1+p})\defi\gamma^{\prime\sharp}(v_0)=\Parens1{\gamma^{\prime\sharp}(v^1_0),\dotsc,\gamma^{\prime\sharp}(v^{1+p}_0)},
  \]
  we obtain a new coordinate system $y\defi(v,\eta)$ on $\aff^{1+p|2q}|_{U_0}$. Notice that
  \[
    v^1=\gamma^{\prime\sharp}(R_0)=\frac R{\gamma^\sharp(R_0)}\gamma^\sharp(R_0)=R.
  \]
  Applying $\gamma^{\prime\sharp}$ to Equation \eqref{eq:polar}, we thus find 
  \begin{equation}\label{eq:super-polar}
    \begin{split}
      u^j&=
      \begin{cases}
        R\sqrt{1-\eta^*\eta}\,\sin(v^{1+p}), &\text{for }j=1+p,\\
        R\sqrt{1-\eta^*\eta}\,\sin(v^j)\cos(v^{j+1})\dotsm\cos(v^{1+p}), &\text{for }j=2,\dotsc,p,\\
        R\sqrt{1-\eta^*\eta}\,\cos(v^2)\cos(v^3)\dotsm\cos(v^{1+p}), &\text{for }j=1,
      \end{cases}\\
      \xi^k&=R\,\eta^k.
    \end{split}
  \end{equation}
  We call these coordinates \Define{super-polar coordinates}. 

  Fixing $R=1$ and noting that $v^1=R$, we obtain a coordinate system on $B|_{B_0\cap U_0}$ by restricting $\tilde y=(\tilde v,\eta)=(v^2,\dotsc,v^{1+p},\eta)$. We shall denote the restrictions by the same symbols. The invariant Berezinian density on $B$ takes the following form. 

  \begin{Lem}[ber-polar]
    On $B_0\cap U_0$, the Berezinian density $\Abs0{Db}$ takes the form 
    \begin{equation}\label{eq:ber-polar}
      \Abs0{Db}=\Abs0{D(\tilde y)}\,\cos(v^3)\cos^2(v^4)\dotsm\cos^{p-1}(v^{1+p})(1-\eta^*\eta)^{\frac{p-1}2}.
    \end{equation}
  \end{Lem}

  \begin{proof}
    A straightforward computation shows
    \[
      \ABer2{\frac{\partial x}{\partial y}}=R^{2\vrho}\cos(v^3)\cos^2(v^4)\dotsm\cos^{p-1}(v^{1+p})(1-\eta^*\eta)^{\frac{p-1}2},
    \]
    so 
    \begin{equation}\label{eq:leb-polar}
      \Abs0{D(x)}=\Abs0{D(y)}\,R^{2\vrho}\cos(v^3)\cos^2(v^4)\dotsm\cos^{p-1}(v^{1+p})(1-\eta^*\eta)^{\frac{p-1}2}
    \end{equation}
    on $U_0$. From this, the shape of $\Abs0{Db}$ follows by applying \cite[Corollary 4.3]{groeger}, since $B=R^{-1}(1)$ and $\frac\partial{\partial R}(R)=1$ with $\frac\partial{\partial R}=\frac{\partial}{\partial v^1}$ defined with respect to $y=(v,\eta)$.
  \end{proof}

  Let $\phi:(0,\infty)\times B\longrightarrow\aff^{1+p|2q}\setminus\{0\}$ be the isomorphism defined by
  \[
    \phi^\sharp(x)\defi r\,x|_B,
  \]
  where $r$ denotes the standard coordinate on $(0,\infty)$. This implies $\phi^\sharp(R)=r$. 

  Define a retraction $\gamma''$ of $(0,\infty)\times B$ by 
  \[
    \phi_0\circ\gamma''=\gamma\circ\phi.
  \]
  Then we compute
  \[
    \gamma''^\sharp(r_0)=r\Norm0u|_B,\quad\gamma''^\sharp(u^j_0|_{B_0})=r\Parens1{\Norm0u^{-1}u^j}|_B.
  \]
  Therefore, by Equations \eqref{eq:leb-polar} and \eqref{eq:ber-polar}, we have 
  \begin{equation}
    \int_{\aff^{1+p|2q}}\Abs0{D(x)}\,f=\int_{(0,\infty)\times B}^{\gamma''}dr\,\Abs0{Db}\,r^{2\vrho}\phi^\sharp(f)(r,b),
  \end{equation}
  in view of \cite[Corollary 2.19]{ap-sph}. We wish to apply the boundary term formula \cite[Theorem 4.5]{ahp-int} and change retractions from $\gamma''$ to ${\id}\times\gamma'|_B$. Since $\gamma''$ does not extend to a retraction on $\aff^1\times B$, we need to shift the boundary, to obtain:
  \begin{equation}\label{eq:polar-bdterms}
    \begin{split}
      \int_{\aff^{1+p|2q}}\Abs0{D(x)}\,f=&\lim_{\eps\to0+}\Biggl[\int_\eps^\infty dr\,r^{2\vrho}\int_B\Abs0{Db}\,\phi^\sharp(f)(r,b)-B_\vrho(f)(\eps)\Biggr]\\ 
      B_\vrho(f)(\eps)\defi&\sum_{k=1}^\infty\frac{\partial^{k-1}_{r=\eps}}{k!}\Parens3{r^{2\vrho+k}\int_B\Abs0{Db}\,(1-\Norm0{b_\ev})^k\phi^\sharp(f)(r,b)},
    \end{split}
  \end{equation}
  where we write $b_\ev=u(b)=(u_1(b),\dotsc,u_p(b))$, $b$ denoting the identity of $B$. Here, the integrals
  \[
    \int_B\Abs0{Db}\,(1-\Norm0{b_\ev})^k\,\phi^\sharp(f)(r,b)
  \]
  where $r$ and $b$ respectively denote the identity of $(0,\infty)$ and $B$, are well-defined functions on $(0,\infty)$, and the sum in the definition of $B_\vrho$ may actually be taken to extend only to $k=q$.

  Notice that if $\vrho\sge0$, then the derivatives in Equation \eqref{eq:polar-bdterms} all converge to zero for $\eps\to0+$. This proves the following statement. 

  \begin{Prop}[ber-polar-i]
    Assume that $\vrho\sge0$. Then for any $f\in\Gamma_c(\sh O_{\smash{\aff^{1+p|2q}}})$, we have
    \begin{equation}\label{eq:ber-polar-i}
      \int_{\aff^{1+p|2q}}^\gamma\Abs0{D(x)}f=\int_0^\infty dr\,r^{2\vrho}\int_B\Abs0{Db}\,\phi^\sharp(f)(r,b).
    \end{equation}
  \end{Prop}

  This gives the desired expression for the invariant Berezin integral on $G/K=D$ in case $\vrho\sge0$.

  \begin{proof}[\prfof{Th}{inv-ber-gk} \eqref{item:inv-ber-gk-i}]
    The assertion follows from Equations \eqref{eq:kak-decomp}, \eqref{eq:aplus-param}, \eqref{eq:km-ber}, \eqref{eq:gk-ber}, and \eqref{eq:ber-polar-i}.
  \end{proof}

  To determine the boundary terms in Equation \eqref{eq:polar-bdterms} more precisely in case and $\vrho<0$, more work has to be done. We Taylor expand 
  \[
    g\defi f-\sum_{\ell=0}^{-2\vrho-1}f_\ell,\quad f_\ell\defi\sum_{\Abs0\alpha=\ell}\frac{x^\alpha}{\alpha!}\partial_{x=0}^\alpha(f)
  \]
  where $\alpha\in\nats^{1+p}\times\{0,1\}^{2q}$ and 
  \[
    \partial^\alpha_x\defi\Parens3{\frac\partial{\partial\xi^{2q}}}^{\alpha_{1+p+2q}}\dotsm\Parens3{\frac\partial{\partial\xi^1}}^{\alpha_{2+p}}\Parens3{\frac\partial{\partial u^{1+p}}}^{\alpha_{1+p}}\dotsm\Parens3{\frac\partial{\partial u^1}}^{\alpha_1}.
  \]
  The homogeneous Taylor components $f_\ell$ for $\Abs0\ell$ odd are by definition anti-invariant under $x\longmapsto-x$, so do not contribute to the boundary terms. Moreover, the Taylor remainder term $\phi^\sharp(g)$ is $r^{-2\vrho}h$ for some superfunction $h$ which extends to $\aff^1\times B$, so that $g$ does also not contribute. We find 
  \begin{equation}\label{eq:polar-bdterms-rhoneg}
      B_\vrho(f)(\eps)=\sum_{k=1}^\infty\sum_{\ell=0}^{\lfloor-\frac12-\vrho\rfloor}\frac{\partial^{k-1}_{r=\eps}(r^{2\vrho+k+2\ell})}{k!}\int_B\Abs0{Db}\,(1-\Norm0{b_\ev})^kf_{2\ell}(b),
  \end{equation}
  by the homogeneity of $f_{2\ell}$. If $\ell<-\frac12-\vrho$, then we have 
  \[
    \frac{\partial_{r=\eps}^{k-1}}{k!}(r^{2\vrho+k+2\ell})=\frac{(-1)^k}{2\vrho+1+2\ell}\binom{-2\vrho-2\ell-1}k\eps^{2\vrho+1+2\ell},
  \]
  and for $\ell=-\frac12-\vrho$, we similarly have $\frac{\partial_{r=\eps}^{k-1}}{k!}(r^{k-1})=\frac1k$. Thus, we obtain
  \[
    \sum_{k=1}^\infty\frac{\partial^{k-1}_{r=\eps}(r^{2\vrho+k+2\ell})}{k!}(1-\Norm0{b_\ev})^k=\eps^{2\vrho+1+2\ell}\cdot
    \begin{cases}
      \displaystyle\frac{\Norm0{b_\ev}^{-2\vrho-2\ell-1}-1}{2\vrho+1+2\ell}, &\text{for }\ell<-\tfrac12-\vrho,\\
      \displaystyle\log(\Norm0{b_\ev}),&\text{for }\ell=-\tfrac12-\vrho.
    \end{cases}
  \]
  Now recall from \thmref{Lem}{ber-polar} that $\Abs0{Db}$ is of the form $\Abs0{D(\tilde y)}\,\gamma'^\sharp(\chi)(1-\eta^*\eta)^{\frac{p-1}2}$ for some smooth function $\chi$ on $B_0$. Since $\xi|_B=\eta|_B$ and $\Abs0{b_\ev}=\sqrt{1-\eta^*\eta}$, the integral over $B$ of 
  \[
    \Abs0{Db}\,\Norm0{b_\ev}^{-2\vrho-2\ell-1}f_{2\ell}(b)=\Abs0{D(\tilde y)}\,\gamma'^\sharp(\chi)(1-\eta^*\eta)^{q-\ell-1}f_{2\ell}(b)
  \]
  vanishes, since the order of the integrand in powers of $\eta$ is at most $2q-2$. 

  For $\vrho<0$ integral, we thus arrive at the equation 
  \begin{equation}\label{eq:polar-bdterms-rhoneg-int}
    B_\vrho(f)(\eps)=-\sum_{\ell=0}^{-1-\vrho}\frac{\eps^{2\vrho+1+2\ell}}{2\vrho+1+2\ell}\int_B\Abs0{Db}\,f_{2\ell}(b).
  \end{equation}
  If $\vrho<0$ is half-integral, then $p$ is odd, so that the integral over $B$ of
  \[
    \Abs0{Db}\,f_{2\ell}(b)=\Abs0{D(\tilde y)}\,(1-\eta^*\eta)^{\frac{p-1}2}f_{2\ell}(b)
  \]
  also vanishes, seeing that $(1-\eta^*\eta)^{\frac{p-1}2}$ has order at most $p-1$ in $\eta$, and that 
  \[
    p-1+2\ell\sle p-2-2\vrho=2q-2<2q.
  \]
  Thus, in this case, we obtain the following expression for the boundary terms:
  \begin{equation}\label{eq:polar-bdterms-rhoneg-halfint}
    B_\vrho(f)(\eps)=\int_B\Abs0{Db}\,\log(\Norm0{b_\ev})f_{-1-2\vrho}(b).
  \end{equation}

  We obtain the following integration formul\ae{}.

  \begin{Prop}[ber-polar-iii]
    Let $\vrho<0$ be half-integral. \Fa $f\in\Gamma_c(\sh O_{\smash{\aff^{1+p|2q}}})$, 
    \begin{equation}\label{eq:ber-polar-iii}
      \begin{split}
        \int_{\aff^{1+p|2q}}\Abs0{D(x)}\,f=\int_0^\infty dr\,&r^{2\vrho}\int_B\Abs0{Db}\,\phi^\sharp(f)(r,b)\\
        &+\frac{\partial^{-1-2\vrho}_{r=0}}{(-1-2\vrho)!}\int_B\Abs0{Db}\,\log(\Norm0{b_\ev})\phi^\sharp(f)(r,b).
      \end{split}
    \end{equation}
  \end{Prop}

  \begin{Prop}[ber-polar-ii]
    Let $\vrho<0$ be integral. \Fa $f\in\Gamma_c(\sh O_{\smash{\aff^{1+p|2q}}})$, 
    \begin{equation}\label{eq:ber-polar-ii}
      \int_{\aff^{1+p|2q}}\Abs0{D(x)}\,f=\frac{(-1)^\vrho\Gamma\Parens1{\vrho+\tfrac12}}{2\sqrt\pi}\int_0^\infty dr\,r^{-\frac12}\partial_r^{-\vrho}\Bracks3{\int_B\Abs0{Db}\,\phi^\sharp(f)(\sqrt{r},b)}.
    \end{equation}
  \end{Prop}

  \begin{Rem}
    The integrals on the right-hand side of Equation \eqref{eq:ber-polar-iii} exist only as iterated integrals in the order written. 
  \end{Rem}

  \begin{proof}[\prfof{Prop}{ber-polar-ii}]
    We may replace $f_{2\ell}$ by $\frac1{(2\ell)!}\partial^{2\ell}_{r=0}\phi^\sharp(f)(r,b)$ in Equation \eqref{eq:polar-bdterms-rhoneg-int}. Moreover, the function $f^\circ$ on $\reals$ given by 
    \[
      f^\circ(r)\defi\frac1{\int_B\Abs0{Db}}\int_B\Abs0{Db}\,\phi^\sharp(f)(r,b)
    \]
    is even, so that $\bar f(r)\defi f^\circ(\sqrt{\Abs0r})$ is a smooth function of $r\in\reals$. Thus, in Equation \eqref{eq:polar-bdterms}, $f$ may be replaced by $h(x)\defi \bar f(\Norm0x^2)$ (\ie the pullback of $\bar f$ under the morphism $\Norm0\cdot^2:\aff^{1+p|2q}\longrightarrow\aff^1$). The claim follows from \cite[Corollary 2.24]{ap-sph}, together with Equation \eqref{eq:volb}.
  \end{proof}

  \begin{proof}[\prfof{Th}{inv-ber-gk} \eqref{item:inv-ber-gk-ii}--\eqref{item:inv-ber-gk-iii}]
    The assertion follows from Equations \eqref{eq:kak-decomp}, \eqref{eq:aplus-param}, \eqref{eq:km-ber}, \eqref{eq:gk-ber}, \eqref{eq:ber-polar-iii}, and \eqref{eq:ber-polar-ii}.
  \end{proof}

  \begin{Cor}[log-int]
    Let $\vrho<0$ be half-integral. Then 
    \begin{equation}\label{eq:log-int}
      \int_{K/M}\Abs0{D\dot k}\,\log\Norm0{ke_1}_\ev=2^q(-\pi)^{\frac{p+1}2}\Gamma\Parens1{\tfrac12-\vrho}c_{K/M}.
    \end{equation}
  \end{Cor}

  \begin{proof}
    Let $\chi\in\Ct[_c^\infty]0{[0,1)}$ be equal to $1$ in a neighbourhood of $0$, and define the compactly supported $K$-invariant superfunction $\tilde\chi\defi\chi(\Norm0x^2)$ on $G/K$. 

    Using Equation \eqref{eq:inv-ber-gk-iii} and the fact that $\vol(K/M)=0$ by Equation \eqref{eq:km-int-norm}, we see
    \begin{align*}
      \int_{K/M}\Abs0{D\dot k}\,\log\Norm0{ke_1}_\ev&=\frac{\partial_{r=0}^{-2\vrho-1}}{(-2\vrho-1)!}\Bracks3{r^{-2\vrho-1}\chi(r^2)\int_{K/M}\Abs0{D\dot k}\,\log\Norm0{ke_1}_\ev}\\
      &=\int_{G/K}\Abs0{D\dot g}\,\Parens1{1-\Norm0{g\cdot 0}^2}^{1+\vrho}\Norm0{g\cdot 0}^{-2\vrho-1}\tilde\chi(\Norm0{g\cdot0}^2).
      \intertext{By Equation \eqref{eq:gk-ber}, this equals}
      &=c_{K/M}\int_{\aff^{1+p|2q}}\Abs0{D(x)}\,\Norm0x^{-2\vrho-1}\tilde\chi(\Norm0x^2).
      \intertext{Applying \cite[Corollary 2.24]{ap-sph}, this equates to}
      &=2^q(-\pi)^{\frac{p+1}2}c_{K/M}\partial_{r=0}^{-\vrho-\frac12}(r^{-\vrho-\frac12}\chi(r))\\
      &=2^q(-\pi)^{\frac{p+1}2}\Gamma\Parens1{\tfrac12-\vrho}c_{K/M},
    \end{align*}
    as was claimed. 
  \end{proof}

  \begin{Rem}
    A direct proof of \thmref{Cor}{log-int} is also not hard. Indeed, we have 
    \[
      \int_{-\pi/2}^{\pi/2}dv^3\cos(v^3)\dotsm\int_{-\pi/2}^{\pi/2}dv^{1+p}\,\cos^{p-1}(v^{1+p})=\frac{\pi^{\frac{p-1}2}}{\Gamma\Parens1{\frac{p+1}2}}
    \]
    by a standard identity for the beta function \cite[Chapter 1.5.1, (19)]{emot}. Moreover, 
    \begin{align*}
      \int_{\aff^{0|2q}}D(\eta)\,&(1-\eta^*\eta)^{\frac{p-1}2}\log(1-\eta^*\eta)=2^qq!\sum_{k+\ell=q,k\sge1}\frac{(-1)^{\ell+1}}k\binom{\frac{p-1}2}\ell.
      \intertext{As $\frac{p-1}2$ is an integer and $q-1\sge\frac{p-1}2$ by the assumption on $\vrho$, this equals}
      &=2^qq!\sum_{\ell=0}^{\frac{p-1}2}\frac{(-1)^{\ell+1}}{q-\ell}\binom{\frac{p-1}2}\ell=(-1)^{\frac{p+1}2}2^q\Gamma\Parens1{\tfrac{p+1}2}\Gamma\Parens1{\tfrac12-\vrho},
    \end{align*}
    in view of the equality
    \[
      \sum_{k=0}^n\binom nk(-1)^k\frac1{q-k}=\frac{(-1)^n}{(q-n)\binom qn}
    \]
    valid for non-negative integers $q>n$. Integrating over $v^2$ contributes a factor of $2\pi$, and the square root in the logarithm of $\Norm0b_\ev$ a factor of $\tfrac12$. Thus, on applying Equation \eqref{eq:ber-polar}, we compute
    \[
      \int_B\Abs0{Db}\log\Norm0b_\ev=2^q(-\pi)^{\frac{p+1}2}\Gamma\Parens1{\tfrac12-\vrho},
    \]
    which plainly gives the desired result. 
  \end{Rem}


  We end this subsection by determining the normalization of $\Abs0{D\dot k}$ in the remaining case of $\vrho<0$ and half-integral. We begin by reducing the problem to one that only depends on $\vrho$, rather than on $p$ and $q$ simultaneously, by localizing the Berezin integral defining $\phi_\lambda^\vrho$. This is possible due to the $M$-invariance of the integrand.

  \begin{Lem}
    Let $c_{p|2q}$ denote the non-zero constant \scth $\Abs0{D\dot k}=c_{p|2q}\Abs0{Db}$ in case $\dim B=p|2q$. Then for $\vrho=\frac p2-q$ negative and half-integral, we have
    \begin{equation}\label{eq:ckmrec}
      c_{p|2q}=\frac1{(-2\pi)^{\frac{p-1}2}}c_{1|2q-p+1}.
    \end{equation}
  \end{Lem}

  \begin{proof}
    Let $j$ be a non-negative integer \scth $2j\sle p=\min(p,2q)$. We may take $j=\frac{p-1}2$. Consider the odd vector fields 
    \[
      Q_n\defi 
      \begin{aligned}[t]
        &u^{1+p-(2n-1)}\frac\partial{\partial\xi^{2q-2(n-1)}}-u^{1+p-2(n-1)}\frac\partial{\partial\xi^{2q-(2n-1)}}\\
        &+\xi^{2q-(2n-1)}\frac\partial{\partial u^{1+p-2(n-1)}}+\xi^{2q-2(n-1)}\frac\partial{\partial u^{1+p-(2n-1)}},
      \end{aligned}
    \]
    $n=1,\dotsc,j$. They commute and are fundamental vector fields for the $M$-action on $\aff^{1+p|2q}$. In particular, they are tangential to $B$. The zero locus of $Q_n$ is seen to equal that of $Q_n^2$, and is given by 
    \[
      u^{1+p-(2n-1)}=u^{1+p-2(n-1)}=0,\quad n=1,\dotsc,j.
    \]
    Fix $r$ and consider 
    \[
      f\defi e^{(\lambda-\vrho)(H(a_r^{-1}k))},
    \]
    $k$ denoting the identity of $K$. Then $f$ is $M$-invariant. The assumptions of \cite[Theorem 3]{schwarz-zaboronsky} are thus seen to be verified for the Berezinian $\Abs0{Db}\,f$. Therefore, there is an even function $\chi$ on $B$, invariant under $Q_1,\dotsc,Q_j$ and equal to $1$ on a neighbourhood of the subsupermanifold $B$ defined by the equations 
    \[
      u^{1+p-(2n-1)}=u^{1+p-2(n-1)}=\xi^{2q-(2n-1)}=\xi^{2q-2(n-1)}=0,\quad n=1,\dotsc,j,
    \]
    and supported away from the north and south pole $\pm e_1\in B_0$. Moreover, for any such function, we have 
    \[
      \int_B\Abs0{Db}\,f=\int_B\Abs0{Db}\,f\chi.
    \]

    Let $y'$ be the super-polar coordinates for $\aff^{1+p-2j|2q-2j}$ and split $x=(u,\xi)$ as $x=(x',x'')$ where $x'\defi(u^1,\dotsc,u^{1+p-2j},\xi^1,\dotsc,\xi^{2q-2j})$. On the open subspace of $\aff^{1+p|2q}|_{\reals^{1+p}\times(-1,1)^{2j}}$ on which $y$ is a system of local coordinates, $(y',x'')$ is a system of local (``cylindrical'') coordinates. 

    \thmref{Lem}{ber-polar} and the considerations before that show that 
    \[
      \begin{split}
        \Abs0{D(x)}&=\Abs0{D(y',x'')}\cos(v^3)\dotsm\cos^{p-1-2j}(v^{1+p-2j})(1-\eta'^*\eta')^{\frac{p-1}2-j}\\
        &=\Abs0{D(R,\tilde y',x'')}\Parens3{\frac{R}{R'}}^{1+p-2j}\prod_{k=3}^{1+p-2j}\cos^{k-2}(v^k)(1-\eta'^*\eta')^{\frac{p-1}2-j}
      \end{split}
    \]
    where $R'=v'^1$ and $\tilde y'=(v'^2,\dotsc,\eta'^{2q-2j})$ if $y'=(v'^1,\dotsc,\eta'^{2q-2j})$. Let $\Abs0{Db'}$ be the Riemannian Berezin density on $B'$. Then
    \[
      \Abs0{Db}=\Abs0{Db'}\Abs0{D(x'')}R'^{-(1+p-2j)}=\Abs0{Db'}\Abs0{D(x'')}(1-H)^{j-\frac{1+p}2}
    \]
    on the open subspace $\dot B$ of $B$ corresponding to $B_0\setminus\{\pm e_1\}$, where
    \[
      H=\sum_{n=1}^jH_n,\quad H_n\defi(u^{1+p-(2n-1)})^2+(u^{1+p-2(n-1)})^2+2\xi^{2q-(2n-1)}\xi^{2q-2(n-1)}.
    \]
    Here, notice that if $\aff''\defi\aff^{2j|2j}|_{(-1,1)^{2j}}$, then $\dot B\cong B'\times\aff''$, the isomorphism being compatible with the splitting of the local coordinate system $(\tilde y',x'')$. Let $\psi$ denote the isomorphism. Then 
    \begin{align*}
      \int_B\Abs0{Db}\,f&=\int_B\Abs0{Db}\,f\chi=\int_{\aff''}\Abs0{D(x'')}\,(1-H)^{j-\frac{1+p}2}\int_{B'}\Abs0{Db'}\,\psi^{-\sharp}(f\chi)(b',x'')\\
      &=(-2\pi)^j\Bracks3{\prod_{n=1}^j\sqrt{1-H_n}\,(1-H)^{j-\frac{1+p}2}\int_{B'}\Abs0{Db'}\,\psi^{-1\sharp}(f\chi)(b',x'')}_{x''=0}\\
      &=(-2\pi)^j\int_{B'}\Abs0{Db'}\,f|_{B'},
    \end{align*}
    where in the second line, we have applied \cite[Lemma 16]{disertori-spencer-zirnbauer} or \cite[Corollary 2.24]{ap-sph}. 

    Recalling the definition of the constants $c_{p|2q}$, we find
    \[
      \phi_\lambda^\vrho(a_r)=c_{p|2q}\int_B\Abs0{Db}\,f=c_{p|2q}(-2\pi)^j\int_{B'}\Abs0{Db'}\,f|_B'=\frac{(-2\pi)^jc_{p|2q}}{c_{p-2j|2q-2j}}\phi_\lambda^\vrho(a_r).
    \]
    Since $r$ was arbitrary, the assertion follows upon setting $j=\frac{p-1}2$.
  \end{proof}

  Having reduced the computation to the case of $p=1$, we now tackle this case. 

  \begin{Lem}
    We have
    \begin{equation}\label{eq:ckmp1}
      c_{1|2q}=\frac{(-1)^q}{\sqrt2\,2\pi}.
    \end{equation}
  \end{Lem}

  \begin{proof}
    First, we compute for $g=kan$, $a=e^{th_0}$, and $\nu^t\defi(1\,e_1^t)$, that 
    \[
      \nu^t\theta(g)^{-1}g\nu=\nu^t\theta(n)^{-1}a^2n\nu=\nu^ta^2\nu=2e^{2t},
    \]
    where the identities $\nu^t\theta(n)^{-1}=\nu^t$ and $n\nu=\nu$ follow from \cite[Equations (4.7-8)]{ap-sph}. Thus, we have 
    \[
      H(g)=\frac12\log\Bracks3{\frac12\nu^t\theta(g)^{-1}g\nu}.
    \]
    If $r=\tanh(t)$, $k$ is the identity of $K$, and $b=ke_1$, then 
    \begin{align*}
      \nu^t\theta(a_r^{-1}k)a_r^{-1}k\nu&=(1\,(ke_1)^*)a_r^{-2}
      \begin{Matrix}1
        1\\ ke_1
      \end{Matrix},
      \intertext{where we use $(ke_1)^*$ as a shorthand for $(ke_1)^t\begin{Matrix}01&0\\ 0&J\end{Matrix}$. Decomposing $b$ according to the coordinates $x=(u,\xi)$ as $b^t=(b_1\,\tilde b)$, we obtain}
      &=2\Parens1{\cosh(t)-\sinh(t)b_1}^2=2\frac{(1-rb_1)^2}{1-r^2},
    \end{align*}
    as we see by applying $\Norm0b^2=b_1^2+\tilde b^*\tilde b=1$, the addition theorem for the hyperbolic functions, and standard expressions for $\tanh^{-1}(r)$. We conclude that 
    \begin{equation}\label{eq:iwasawa-h}
      H(a_r^{-1}k)=\log\Parens1{\cosh(t)-\sinh(t)b_1}=\log\Bracks3{\frac{1-rb_1}{\sqrt{1-r^2}}}.      
    \end{equation}

    In order to evaluate $c_{1|2q}$, we let $r=\tanh(t)$ and define
    \[
      f^q_\lambda(t)\defi\int_B\Abs0{Db}\,e^{(\lambda-\frac12+q)(H(a_r^{-1}k))},
    \]
    so that $\phi_\lambda^{\frac12-q}(t)=c_{1|2q}f^q_\lambda(t)$. In view of Equations \eqref{eq:ber-polar} and \eqref{eq:super-polar}, and by the above considerations, we have 
    \[
      e^{-(\lambda-\vrho)t}f^q_\lambda(r)=\int_{-\pi}^\pi dv\int D(\eta)\,\Parens1{\tfrac12(1+e^{-2t}-(1-e^{-2t})\cos(v)\sqrt{1-\eta^*\eta})}^{\lambda-\frac12+q}.
    \]
    This entails
    \[
      \frac{\mathbf c_{\frac12-q}(\lambda)}{c_{1|2q}}=2^{\frac32-\lambda-q}\int_0^\pi dv\int_{\aff^{0|2q}}D(\eta)\,\Parens1{1-\cos(v)\sqrt{1-\eta^*\eta}}^{\lambda-\frac12+q}=I_1+I_2,
    \]    
    where we set 
    \[
      \begin{split}
        I_1&\defi 2^{\frac32-\lambda-q}\int_0^{\frac\pi2} dv\int_{\aff^{0|2q}}D(\eta)\,\Parens1{1-\cos(v)\sqrt{1-\eta^*\eta}}^{\lambda-\frac12+q},\\
        I_2&\defi 2^{\frac32-\lambda-q}\int_0^{\frac\pi2} dv\int_{\aff^{0|2q}}D(\eta)\,\Parens1{1+\sin(v)\sqrt{1-\eta^*\eta}}^{\lambda-\frac12+q}.
      \end{split}
    \]
    To evaluate $I_1$, using Equation \eqref{eq:super-polar}, we transform to $u=\cos(v)\sqrt{1-\eta^*\eta}$ and $\xi=\eta$. Since
    \[
      \ABer3{\frac{\partial(u,\xi)}{\partial(v,\eta)}}=\Abs0{\sin(v)}\sqrt{1-\eta^*\eta}=\sqrt{1-u^2-\xi^*\xi},
    \]
    we obtain
    \[
      I_1=2^{\frac32-\lambda-q}\int_0^1du\,(1-u)^{\lambda-\frac12+q}\int D(\xi)\,(1-u^2-\xi^*\xi)^{-\frac12}.      
    \]
    If $k$ is an integer, then decomposing $\xi=(\xi^1,\xi^2,\xi')$, we find
    \[
      \begin{aligned}
        \int D(\xi)\,(1-u^2-\xi^*\xi)^{\frac k2}&=-k\int D(\xi^2,\xi')\,(1-u^2-\xi'^*\xi')^{\frac k2-1}\xi^2\cr
        &=-k\int D(\xi')\,(1-u^2-\xi'^*\xi')^{\frac k2-1},            
      \end{aligned}
    \]
    so that the fermionic integral equals
    \[
      \int D(\xi)\,(1-u^2-\xi^*\xi)^{-\frac12}=(-2)^qq!\binom{-\frac12}q(1-u^2)^{-\frac12-q}=\frac{(-2)^q\sqrt\pi}{\Gamma\Parens1{\frac12-q}}(1-u^2)^{-\frac12-q}.
    \]
    Inserting this back in our expression, we find
    \[
      I_1=\frac{(-1)^q2^{\frac32-\lambda}\sqrt\pi}{\Gamma\Parens1{\frac12-q}}\int_0^1du\,(1-u)^{\lambda-1}(1+u)^{-\frac12-q}.
    \]
    Using the substitution $u=\sin(v)\sqrt{1-\eta^*\eta}$, $\xi=\eta$, a similar computation shows
    \[
      I_2=\frac{(-1)^q2^{\frac32-\lambda}\sqrt\pi}{\Gamma\Parens1{\frac12-q}}\int_0^1du\,(1+u)^{\lambda-1}(1-u)^{-\frac12-q}.
    \]
    A standard identity for the beta function \cite[Chapter 1.5.1, (10)]{emot} yields
    \[
      \frac{\mathbf c_{\frac12-q}(\lambda)}{c_{1|2q}}=\frac{(-1)^q2^{1-q}\sqrt\pi}{\Gamma\Parens1{\frac12-q}}B\Parens1{\lambda,\tfrac12-q}=\frac{(-1)^q2^{1-q}\sqrt\pi\,\Gamma(\lambda)}{\Gamma\Parens1{\lambda+\frac12-q}}.
    \]
    Comparing this with Equation \eqref{eq:cfn-gkfmla} gives the desired result. 
  \end{proof}

  We have proved the following proposition. 

  \begin{Prop}[norm-const]
    Let $\vrho<0$ be half-integral. Then 
    \begin{equation}\label{eq:norm-const}
      c_{K/M}=\frac{(-1)^q}{\sqrt2(2\pi)^{\frac{p+1}2}},    
    \end{equation}
    and in particular,
    \begin{equation}\label{eq:log-int-2}
      \int_{K/M}\Abs0{D\dot k}\,\log\Norm0{ke_1}_\ev=(-1)^{\vrho+\frac12}2^{-\vrho-1}\Gamma\Parens1{\tfrac12-\vrho}.
    \end{equation}
  \end{Prop}

  \begin{proof}
    The latter equation follows from the former by inserting the value of $c_{K/M}$ into Equation \eqref{eq:log-int}, and the former follows from Equations \eqref{eq:ckmrec} and \eqref{eq:ckmp1}.
  \end{proof}

  \subsection{Non-Euclidean Fourier inversion}

  In this subsection, we define the non-Euclidean Fourier (or Helgason--Fourier) transform for functions on the super-hyperbolic space $D=G/K$ and prove an inversion formula for this transform. 

  \begin{Def}[fourier][Helgason--Fourier transform]
    Let $f\in\Gamma_c(\sh O_{G/K})=\Gamma_c(\sh O_G)^K$. We define $\sh F(f)\in\Gamma(\sh O_{\aff(i\ger a_\reals^*)\times K/M})=\Gamma(\sh O_{\aff(i\ger a_\reals^*)\times K})^M$ by 
    \begin{equation}\label{eq:fourier}
      \sh F(f)(\lambda,k)\defi\int_{G/K}\Abs0{D\dot g}\,f(g)e^{(\lambda-\vrho)(H(g^{-1}k))},
    \end{equation}
    where $\lambda$ and $k$ denote the identity of $\aff(i\ger a_\reals^*)$ and $K$, respectively. The map 
    \[
      \sh F:\Gamma_c(\sh O_{G/K})\longrightarrow\Gamma(\sh O_{\aff(i\ger a_\reals^*)\times K/M})  
    \]
    is called the \Define{Helgason--Fourier transform} or \Define{non-Euclidean Fourier transform}. Note that as in the definition of the wave packet transform, we are using the Iwasawa $A$ projection $H(-)$ corresponding to the $KAN$ decomposition, rather than using the projection $A(-)$ associated with the $NAK$ decomposition, as Helgason does. 
  \end{Def}

  \begin{Rem}
    By \thmref{Lem}{gk-km-ber} and by using a workable expression for $H$, for instance the one given in \cite[Lemma 4.9]{ap-sph}, the non-Euclidean Fourier transform can be written in totally explicit form. This can also be achieved by embedding $D=G/K$ equivariantly as a hypersurface (a one-sheeted hyperboloid) in $\aff^{2+p|2q}$.
  \end{Rem}

  In order to formulate our main theorem succinctly, recall the \Define{convolution} 
  \begin{equation}\label{eq:convolution-def}
    (f_1*f_2)(h)\defi\int_{G/K}\Abs0{D\dot g}f_1(hg)f_2(g^{-1})=\int_{G/K}\Abs0{D\dot g}\,f_1(g)f_2(g^{-1}h),
  \end{equation}
  defined \fa supermanifolds $T$ and all $g\in_TG$ if $f_1,f_2\in\Gamma(\sh O_{G/K})$ where at least one is compactly supported and $f_2$ is $K$-biinvariant. Moreover, recall the wave packet transform $\sh J$ from Equation \eqref{eq:wavepacket-def}.

  \begin{Th}[fouinv][Fourier inversion formula]
    Let $f\in\Gamma_c(\sh O_{G/K})=\Gamma_c(\sh O_G)^K$. Then $\sh J(\sh F(f))$ exists and we have the following inversion formul\ae.
    \begin{enumerate}[wide]
      \item\label{item:fouinv-i} Let $\vrho\sge0$. Then we have 
      \begin{equation}\label{eq:fouinv-i}
        2^{2(1-\vrho)}\pi f=\sh J(\sh F(f)).
      \end{equation}
      \item\label{item:fouinv-ii} Let $\vrho<0$. Then we have 
      \begin{equation}\label{eq:fouinv-ii}
        2^{2(1-\vrho)}\pi\,f=\sh J(\sh F(f))-f*\sh J(1).
      \end{equation}
    \end{enumerate}
  \end{Th}

  Here, notice that we take the liberty to consider functions on $G/K$ either as $K$-invariant functions on $G$ or as functions on $D$.

  \begin{Rem}
    Observe that \thmref{Th}{fouinv} recovers both Helgason's classical inversion formula \cite[Chapter III, Theorem 1.3]{helgason-sym} in the case of the hyperbolic space (\ie $q=0$, and in particular $\vrho>0$), and also Zirnbauer's inversion formula for the super-Poincar\'e disc \cite[Corollary to Theorem 2]{zirn-cmp}. Indeed, in the latter case, we have $\vrho=-\tfrac12$ and $\sh J(1)$ is just a constant, by \thmref{Cor}{wavepacket-const}. To compare with Helgason's formula, note the shift by $2\vrho$ introduced by switching between the Iwasawa $A$ projections corresponding to the $KAN$ and $NAK$ decompositions.
  \end{Rem}

  \begin{proof}[\prfof{Th}{fouinv}]
    Recall from \thmref{Cor}{wavepacket-const} that $\sh J(1)$ is well-defined for $\vrho<0$. To unify the different cases as much as possible, we define $\sh J(1)\defi 0$ for $\vrho\sge0$. As a first step, observe that it is sufficient to prove the identity
    \begin{equation}\label{eq:fouinv-origin}
      2^{2(1-\vrho)}\pi\,f(0)=\sh J(\sh F(f))(0)-(f*\sh J(1))(0)    
    \end{equation}
    \fa functions $f$. Indeed, assume that such an equation has been shown. Then we may introduce an auxiliary supermanifold $T$ and observe that the fibrewise Berezin integral is continuous on the space of fibrewise compactly supported sections of the sheaf $\Abs0{\sh Ber}_{(T\times G/K)/T}$ with its standard LF topology, compare \cite[2.6]{ap-sph} and \cite[Appendix C]{as-sbos} for the relevant definitions. 

    Since $\Gamma(\sh O_T)\otimes\Gamma_c(\Abs0{\sh Ber}_{G/K})$ is dense in this space (compare \cite[Corollary C.9]{as-sbos}), it follows that Equation \eqref{eq:fouinv-origin} holds for the corresponding parameter versions of $\sh F$ and $\sh J$ and all fibrewise compactly supported sections $f$ of $\sh O_{T\times G/K}$.

    Fix $f\in\Gamma_c(\sh O_{G/K})$ and $h\in_TG$. Define $f^h\in\Gamma(T\times G/K)$ by $f^h(g)\defi f(hg)$, \ie 
    \[
      f^h=(h\times{\id}_{G/K})^\sharp(m^\sharp(f)),
    \]
    where $m:G\times G/K\longrightarrow G/K$ denotes the action of $G$ on $G/K$. Then $f^h$ is fibrewise compactly supported, and Equation \eqref{eq:fouinv-origin} applies by assumption. 

    On the other hand, we have 
    \[
      \phi_\lambda(g^{-1}h)=\int_{K/M}\Abs0{D\dot k}\,e^{(\lambda-\vrho)(H(g^{-1}k))}e^{-(\lambda+\vrho)(H(h^{-1}k))}
    \]
    by \cite[Proposition 3.19]{ap-sph}, where $g$ denotes the identity of $G$. Equation \eqref{eq:wavepacket-def} implies
    \[
      \begin{split}
        \sh J(\sh F(f))(h\cdot 0)&=\int_{i\ger a_\reals^*}\frac{d\lambda}{\Abs0{\mathbf c(\lambda)}^2}\int_{G/K}\Abs0{D\dot g}\,f(g)\phi_\lambda(g^{-1}h)\\
        &=\int_{i\ger a_\reals^*}\frac{d\lambda}{\Abs0{\mathbf c(\lambda)}^2}\int_{G/K}\Abs0{D\dot g}\,f^h(g)\phi_\lambda(g^{-1})=\sh J(\sh F(f^h))(0).
      \end{split}
    \]
    Here, we have used the left-invariance of $\Abs0{D\dot g}$ \cite[Lemma 3.11]{ap-sph}. Similarly, we find that $(f^h*\sh J(1))(0)=(f*\sh J(1))(h\cdot 0)$. Therefore, the inversion formula at an arbitrary $h\in_TG$ (for arbitary $T$) follows already from Equation \eqref{eq:fouinv-origin}. In particular, for $T=G$ and $h={\id}_G$, we obtain the statement of \thmref{Th}{fouinv}.

    \bigskip\noindent
    Next, we consider the case of $\vrho\sge0$ or $\vrho<0$ and $\vrho$ integral. To that end, we introduce the abbreviations 
    \begin{equation}
      \psi_\lambda^\vrho(r)\defi\frac1{\mathbf c(\lambda)\mathbf c(-\lambda)}\phi_\lambda^\vrho(a_{\sqrt r}),\quad\Psi^\vrho(r)\defi\sh J(1)(a_{\sqrt r})=\int_{i\ger a_\reals}d\lambda\,\psi^\vrho_\lambda(\sqrt r).
    \end{equation}
    Here, recall the definition of $a_r$ from Equation \eqref{eq:aplus-param}. We will now briefly suspend the proof of the theorem and establish the following two propositions for a general compactly supported smooth function $h$ on $[0,1)$.  
  \end{proof}

  In the following two propositions, notice that both sides of the equations depend only on $\vrho$ and not on $p$ and $q$ individually. This will enable us to prove them by induction on $\vrho$.

  \begin{Prop}[fouinv-origin-i]
    Let $\vrho\sge0$. Then 
    \begin{equation}\label{eq:fouinv-origin-i}
        2^{2(1-\vrho)}\pi\,h(0)=\frac{2^{-\vrho-1}}{\Gamma\Parens1{\vrho+\frac12}}\int_{i\ger a^*}d\lambda\int_0^1dr\,r^{\vrho-\frac12}\psi_\lambda^\vrho(r)h(r),
    \end{equation}
    for any compactly supported smooth function $h$ on $[0,1)$.
  \end{Prop}

  Observe that the above constant $2^{-\vrho-1}\Gamma\Parens1{\tfrac12+\vrho}$ is precisely $\tfrac12\vol(K/M)$.

  \begin{Prop}[fouinv-origin-ii]
    Let $\vrho<0$ be integral. Then 
    \begin{equation}\label{eq:fouinv-origin-ii}
        (-2)^{-\vrho}8\pi^{\frac32}\,h(0)+\int_0^1\!\frac{dr}{\sqrt r}\,\partial^{-\vrho}_r\Parens1{(\Psi^\vrho h)(r)}=\int_{i\ger a^*}\!d\lambda\int_0^1\!\frac{dr}{\sqrt r}\,\partial_r^{-\vrho}\Parens1{(\psi_\lambda^\vrho h)(r)}
    \end{equation}
    for any compactly supported smooth function $h$ on $[0,1)$.
  \end{Prop}

  In the \emph{proof} of these propositions, we need the following lemma. 

  \begin{Lem}[psi-der]
    Let $\vrho$ be arbitrary. Then for $r\in[0,1)$, we have 
    \begin{equation}\label{eq:psi-der}
      \partial_r\psi_\lambda^\vrho(r)=-2(1-r)^{-\frac32}\psi_\lambda^{\vrho+1}(r).
    \end{equation}
    If in addition $\vrho<-1$, then for $r\in[0,1)$, we have 
    \begin{equation}\label{eq:Psi-der}
      \partial_r\Psi^\vrho(r)=-2(1-r)^{-\frac32}\Psi^{\vrho+1}(r).
    \end{equation}
  \end{Lem}

  \begin{proof}
    For $r=\tanh^2(t)\sge0$, we have
    \[
      \cosh(t)=\frac1{\sqrt{1-r}},\quad\sinh(t)=\sqrt{\frac{r}{1-r}},
    \]
    so that 
    \[
      \slashed\partial h(\tanh^2(t))=2(1-r)^{\frac32}\partial_rh(r)
    \]
    for any differentiable function $h(r)$. Applying this in Equation \eqref{eq:sph-der} establishes Equation \eqref{eq:psi-der}. This does not immediately imply Equation \eqref{eq:Psi-der}, since one cannot in general exchange integral and derivative. 

    However, we may apply the second identity in Equation \eqref{eq:wavepacket-const}. Setting 
    \[
      \chi_\vrho(t)\defi\frac{\partial^{-2\vrho-1}_{s=0}}{\Gamma(-2\vrho)}\Bracks3{\frac{(1-2s^2\cosh(t)+s^4)^{-\vrho}}{(1-s)^2}},
    \]
    we observe
    \[
      (\slashed\partial\chi_\vrho)(t)=2\vrho\,\chi_{\vrho+1}(t).
    \]
    This immediately gives Equation \eqref{eq:Psi-der}. Alternatively, one may use the first identity in Equation \eqref{eq:wavepacket-const}. In any case, this proves the claim. 
  \end{proof}

  \begin{proof}[\prfof{Prop}{fouinv-origin-i}]
    Let $\vrho\sge1$. Applying Equation \eqref{eq:psi-der}, integration by parts shows that 
    \begin{align*}
      \int_0^1dr\,r^{\vrho-\frac12}\psi_\lambda^\vrho(r)h(r)
      &=-\frac12\int_0^1\frac{dr\,r^{\vrho-\frac12}}{(1-r)^{\smash{\frac32}}}\partial_r\psi^{\vrho-1}_\lambda(r)h(r)\\
      &=\frac12\int_0^1dr\,r^{\vrho-\frac32}\psi^{\vrho-1}_\lambda(r)\tilde h(r),
    \end{align*}
    where 
    \[
      \tilde h(r)\defi\Parens1{\vrho-\tfrac12}\frac{h(r)}{(1-r)^{\smash{\frac32}}}+r\,\partial_r\Bracks3{\frac{h(r)}{(1-r)^{\smash{\frac32}}}}.
    \]
    Assume that the statement has been proved for $\vrho-1$. Then 
    \[
      \begin{split}
        \frac{2^{-\vrho-1}}{\Gamma\Parens1{\frac12+\vrho}}\int_{i\ger a^*}\!d\lambda&\int_0^1dr\,r^{\vrho-\frac12}\psi^\vrho_\lambda(r)h(r)\\
        &=\frac{2^{-(\vrho-1)}}{4\Gamma\Parens1{\tfrac12+\vrho}}\int_{i\ger a^*}\!d\lambda\int_0^1dr\,r^{\vrho-\frac32}\psi^{\vrho-1}_\lambda(r)\tilde h(r)\\
        &=\frac{2^{2(1-(\vrho-1))}\pi}{4\Parens1{\vrho-\frac12}}\tilde h(0)=2^{2(1-\vrho)}\pi h(0).
      \end{split}
    \]
    This reduces the proof the cases of $\vrho=0$ and $\vrho=\tfrac12$.

    Let $\vphi(r)\defi h(r)(1-r)^{1+\vrho}$ and observe that 
    \[
      \int_0^1\!\frac{dr\,r^{\vrho-\frac12}}{(1-r)^{1+\vrho}}\psi_\lambda^\vrho(r)\vphi(r)=2\int_0^\infty\!dt\,\sinh^{2\vrho}(t)\,(\psi^\vrho_\lambda\vphi)(\tanh^2(t)).
    \]
    Beginning with the case of $\vrho=0$, for $r=\tanh^2(t)$, Equations \eqref{eq:sph-012} and \eqref{eq:cfn-012} give 
    \[
      \psi^0_\lambda(r)=4\sqrt\pi\cosh(\lambda t),
    \]
    at least for $\lambda\notin\frac12\ints$. So, with
    \[
      \Phi(s)=\int_{-\infty}^\infty dt\,e^{-ist}\vphi(\tanh^2(t))=2\int_0^\infty dt\,\cosh(ist)\vphi(\tanh^2(t))
    \]
    denoting the Fourier transform of $\vphi(\smash{\tanh^2(t)})$, the inversion formula for the Fourier transform gives 
    \[
      \begin{split}
        \int_{i\ger a^*}\!d\lambda\int_0^1dr\,r^{-\frac12}\psi_\lambda^0(r)h(r)&=4\sqrt\pi\int_{-\infty}^\infty\!ds\,\Phi(s)=8\pi\sqrt\pi\,\vphi(\tanh^2(0))=8\pi\sqrt\pi\,h(0).
      \end{split}
    \]
    As $2^{3-\vrho}\Gamma\Parens1{\vrho+\tfrac12}=8\sqrt\pi$ for $\vrho=0$, this proves the claim in this case.

    In the case of $\vrho=\frac12$, Equations \eqref{eq:sph-012} and \eqref{eq:cfn-012} give, for $\lambda\notin\frac12\ints$,
    \[
      \begin{split}
        \psi^{\frac12}_\lambda(r)&=\pi\sqrt2\,\frac{\Gamma\Parens1{\lambda+\frac12}\Gamma\Parens1{-\lambda+\frac12}}{\Gamma(\lambda)\Gamma(-\lambda)}P_{\smash{\lambda-\frac12}}(\cosh(t))\\
        &=-\pi\sqrt2\,\lambda\tan(\pi\lambda)P_{\smash{\lambda-\frac12}}(\cosh(t)),
      \end{split}
    \]
    where in the last identity, we have applied the classical Euler reflection formula. 

    Thus, we obtain with $\vphi$ as above and the substitution $u=\cosh(t)$, that 
    \[
      \begin{split}
        \int_{i\ger a^*}\!d\lambda &\int_0^1\!dr\,\psi_\lambda^{\frac12}(r)h(r)\\
        &=2\sqrt2\pi\int_{-\infty}^\infty\!ds\,s\tanh(\pi s)\int_1^\infty du\,P_{\smash{-\frac12}+is}(u)\vphi\Parens1{1-\tfrac1{u^2}}\\
        &=4\sqrt2\pi \vphi\Parens2{1-\frac1{u^2}}\Big|_{u=1}=4\sqrt2\pi h(0),
      \end{split}
    \]
    by the inversion formula for the Mehler--Fock transform \cite[Chapter 3.14, (8--9)]{emot}. As $2^{3-\vrho}\Gamma\Parens1{\vrho+\tfrac12}=4\sqrt2$ for $\vrho=\frac12$, this proves the claim in this case.
  \end{proof}

  \begin{Rem}
    \thmref{Prop}{fouinv-origin-i} can be derived from the general inversion formula for the non-Euclidean Fourier transform of $K$-invariant functions on hyperbolic space \cite[Chapter IV, Theorem 7.5]{helgason-gga} once it is clear that the integrals only depend on $\vrho$. However, we prefer to give the above proof, which, in addition to being self-contained, has the merit that its technique also applies in the case of $\vrho<0$. 
  \end{Rem}

  \begin{proof}[\prfof{Prop}{fouinv-origin-ii}]
    First, we claim that for $\vrho\sle-1$, we have
    \begin{equation}\label{eq:psi-Psi}
      \int_{i\ger a^*}\!d\lambda\int_0^1dr\,r^{-\frac12}\psi_\lambda^\vrho(r)h(r)=\int_0^1dr\,r^{-\frac12}\Psi^\vrho(r)h(r).
    \end{equation}
    In view Equation \eqref{eq:sphfn-wt-der-est}, this is certainly the case if $\vrho\sle-2$, since we may exchange the order of integration. To prove the claim for $\vrho=-1$, let $H(r)\defi-\int_r^1ds\,s^{-\frac12}h(s)$. Using integration by parts and Equation \eqref{eq:psi-der}, we find
    \begin{align*}
      \int_{i\ger a^*}\!d\lambda\int_0^1&dr\,r^{-\frac12}\psi_\lambda^{-1}(r)h(r)=\int_{i\ger a^*}\!d\lambda\int_0^1dr\,\psi_\lambda^{-1}(r)H'(r)\\
      &=-\int_{i\ger a^*}\!d\lambda\,\psi_\lambda^{-1}(0)H(0)+2\int_{i\ger a^*}\!d\lambda\int_0^1dr\,(1-r)^{-\frac32}\psi^0_\lambda(r)H(r).
      \intertext{Applying Equation \eqref{eq:fouinv-origin-i}, we find}
      &=-\Psi^{-1}(0)H(0)+16\pi^{\frac32}\Bracks1{r^{\frac12}(1-r)^{-\frac32}H(r)}_{r=0}.
      \intertext{Evaluating Equation \eqref{eq:wavepacket-const} for $\vrho=-1$ yields $\Psi^{-1}(r)\equiv\text{const}.$, so that this equals}
      &=\int_0^1dr\,r^{-\frac12}\Psi^{-1}(r)h(r),
    \end{align*}
    as was claimed. This establishes Equation \eqref{eq:psi-Psi}.

    To prove the proposition, we define 
    \[
      T_\vrho(h)\defi\int_{i\ger a^*}d\lambda\int_0^1dr\,r^{-\frac12}\partial_r^{-\vrho}\Parens1{\psi_\lambda^\vrho(r)h(r)}-\int_0^1dr\,r^{-\frac12}\partial^{-\vrho}_r\Parens1{\Psi^\vrho(r)h(r)}.
    \]
    Assume that $\vrho\sle-2$ and the assertion has been proved for negative integers $>\vrho$. 

    Observe that by Equation \eqref{eq:psi-der}
    \[
      \partial_r^{-\vrho}(\psi_\lambda^\vrho(r)h(r))=-2\partial_r^{-\vrho-1}\Parens1{(1-r)^{-\frac32}\psi_\lambda^{\vrho+1}(r)h(r)}+\sum_{k=0}^{-\vrho-1}\psi_\lambda^{\vrho+k}(r)h_k(r)
    \]
    where $h_k$ is a polynomial in the derivatives of $h$ and $(1-r)^{-\frac12}$, and so is a compactly supported smooth function on $[0,1)$. The same equation holds for $\Psi^\vrho$ in place of $\psi_\lambda^\vrho$, by Equation \eqref{eq:Psi-der}. By Equation \eqref{eq:psi-Psi}, it follows that 
    \[
      T_\vrho(h)=-2\,T_{\vrho+1}\Parens1{(1-r)^{-\frac32}h(r)}=\dotsm=(-2)^{-\vrho-1}T_{-1}\Parens1{(1-r)^{\frac32(\vrho+1)}h(r)},
    \]
    by the inductive assumption. This reduces the assertion to the case of $\vrho=-1$.

    In this case, since $\Psi^{-1}$ is constant, Equations \eqref{eq:psi-der}, \eqref{eq:psi-Psi}, and \eqref{eq:fouinv-origin-i} give
    \[
      T_{-1}(h)=-2\int_{i\ger a^*}\!d\lambda\int_0^1dr\,r^{-\frac12}(1-r)^{-\frac32}\psi_\lambda^0(r)h(r)=-16\pi\sqrt\pi\,h(0).
    \]
    This proves the assertion. 
  \end{proof}

  \begin{proof}[\prfof{Th}{fouinv} (continued)]
    We now define, for $f\in\Gamma_c(\sh O_{G/K})$, 
    \[
      h(r)\defi(1-r)^{-1-\vrho}\int_B\Abs0{D\dot k}\,f(ka_{\sqrt r}).
    \]
    In case $\vrho\sge0$, \thmref{Th}{inv-ber-gk} gives 
    \begin{align*}
      \sh J(\sh F(f))(0)&=\int_{i\ger a_\reals}\frac{d\lambda}{\Abs0{\mathbf c(\lambda)}^2}\int_B\Abs0{D\dot k}\int_{G/K}f(g)e^{(\lambda-\vrho)(H(g^{-1}k))}\\
      &=\frac12\int_{i\ger a_\reals}d\lambda\int_0^1dr\,r^{\vrho-\frac12}\psi_\lambda^\vrho(r)h(r),
      \intertext{which by \thmref{Prop}{fouinv-origin-i} equals}
      &=2^{2-\vrho}\pi\,\Gamma\Parens1{\vrho+\tfrac12}h(0)=2^{2(1-\vrho)}\pi\,f(0),
    \end{align*}
    where in the last step, we have applied Equation \eqref{eq:km-int-norm}.

    When $\vrho<0$ is integral, we similarly compute
    \begin{align*}
      \sh J(\sh F(f))(0)&=\frac{(-1)^\vrho\Gamma\Parens1{\vrho+\frac12}}{2\sqrt\pi}\int_{i\ger a_\reals^*}d\lambda\int_0^1dr\,r^{-\frac12}\partial_r^{-\vrho}\Bracks1{\psi_\lambda^\vrho(r)h(r)},
      \intertext{which by \thmref{Prop}{fouinv-origin-ii} equals}
      &=2^{-\vrho+2}\pi\,\Gamma\Parens1{\vrho+\tfrac12}h(0)+\frac{(-1)^q\Gamma\Parens1{\vrho+\frac12}}{2\sqrt\pi}\int_0^1dr\,r^{-\frac12}\partial_r^{-\vrho}\Bracks1{\Psi^\vrho(r)h(r)}.
      \intertext{Applying Equation \eqref{eq:km-int-norm} and \thmref{Th}{inv-ber-gk} once again, together with the fact that $\sh J(1)$ is $K$-biinvariant and $\sh J(1)(g^{-1})=\sh J(1)(g)$, we find}
      &=2^{2(1-\vrho)}\pi\,f(0)+(f*\sh J(1))(0).
    \end{align*}
    This proves \thmref{Th}{fouinv} in case $\vrho\sge0$ or $\vrho<0$ is integral.
  \end{proof}

  To complete the proof in case $\vrho<0$ is half-integral, we need the following lemma.

  \begin{Lem}
    Let $\vrho<0$ be half-integral. Then
    \begin{align}\label{eq:psi-der-zero}
      \partial^k_{r=0}\psi^\vrho_\lambda(r)=0\qquad\forall k<-\vrho+\tfrac12,\\
      \label{eq:Psi-der-zero}
      \partial^k_{r=0}\Psi^\vrho(r)=0\qquad\forall k<-\vrho-\tfrac12.
    \end{align}
    Moreover, we have 
    \begin{equation}\label{eq:Psi-topder}
      \partial^{-\vrho-\frac12}_{r=0}\Psi^\vrho(r)=(-1)^{\vrho-\frac12}2^{3-\vrho}\pi.
    \end{equation}
  \end{Lem}

  \begin{proof}
    By the definitions, we have $\psi^\vrho_\lambda(0)=0$ for $\vrho\sle-\tfrac12$. In view of Equation \eqref{eq:psi-der}, the derivatives $\partial_{r=0}^k\psi^\vrho_\lambda(r)$ for $1\sle k<-\vrho+\frac12$ for are linear combinations of the derivatives $\partial_{r=0}^\ell\psi^{\vrho+1}_\lambda(r)$ for $0\sle\ell\sle k-1$, \ie for $\ell<-(\vrho+1)+\tfrac12$. Thus, by induction on $\vrho$, the proof of Equation \eqref{eq:psi-der-zero} is reduced to the case of $\vrho=-\tfrac12$, which we have already established.

    Similarly, we have $\Psi^\vrho(0)=0$ for $\vrho\sle-\tfrac32$, by Equation \eqref{eq:wavepacket-const-halfint}. By Equation \eqref{eq:Psi-der}, the derivatives $\partial_{r=0}^k\Psi^\vrho(r)$ for $1\sle k<-\vrho-\frac12$ for are linear combinations of the derivatives $\partial_{r=0}^\ell\Psi^{\vrho+1}(r)$ for $0\sle\ell\sle k-1$, \ie for $\ell<-(\vrho+1)-\tfrac12$. Thus, by induction on $\vrho$, the proof of Equation \eqref{eq:Psi-der-zero} is reduced to the case of $\vrho=-\tfrac12$. But in this case, there is nothing to prove, since the assumption on $k$ is never fulfilled. 

    For the final equation, notice that $r=\tanh^2(t)$ gives $1-\cosh(t)=1-(1-r)^{-\frac12}$. In view of Equation \eqref{eq:wavepacket-const-halfint} and by the above considerations, 
    \[
      \partial_{r=0}^k\Parens1{1-(1-r)^{-\frac12}}^{-\vrho-\frac12}=0,\quad k=0,\dotsc,-\vrho-\tfrac32.
    \]
    This implies 
    \[
      \begin{split}
        \frac1{\Parens1{-\vrho-\frac12}!}&\partial_{r=0}^{-\vrho-\frac12}\Parens1{1-(1-r)^{-\frac12}}^{-\vrho-\frac12}=\lim_{r\to0}r^{\vrho+\frac12}\Parens1{1-(1-r)^{-\frac12}}^{-\vrho-\frac12}\\
        &=(-1)^{\vrho+\frac12}\Bracks3{\lim_{r\to 0}\frac{(1-r)^{-\frac12}-1}r}^{-\vrho-\frac12}=(-2)^{\vrho+\frac12}.
      \end{split}
    \]
    Inserting this back in Equation \eqref{eq:wavepacket-const-halfint} gives the claim.
  \end{proof}

  \begin{proof}[\prfof{Th}{fouinv} (continued)]
    Let $\vrho<0$ be half-integral. Let $f\in\Gamma_c(\sh O_{G/K})$ be arbitrary, and define 
    \[
      h_1(r)\defi\int_{K/M}\Abs0{D\dot k}\,f(ka_{\sqrt r}),\quad h_2(r)\defi\int_{K/M}\Abs0{D\dot k}\,\log\Norm0{ke_1}_\ev\,f(ka_{\sqrt r}).
    \]
    Then Equation \eqref{eq:inv-ber-gk-iii} may be rephrased as follows: 
    \[
      \int_{G/K}\Abs0{D\dot g}\,f(g)=2\int_0^1dr\,r^{\vrho+1}(1-r)^{-1-\vrho}h_1(r)+\frac{\partial_{r=0}^{-\vrho-\frac12}}{\Gamma\Parens1{\frac12-\vrho}}\Bracks1{(1-r)^{-1-\vrho}h_2(r)}.
    \]
    Applying Equation \eqref{eq:psi-der-zero}, we therefore find 
    \[
      \sh J(\sh F(f))(0)=2\int_{i\ger a^*}\!d\lambda\int_0^1dr\,r^{\vrho+1}(1-r)^{-1-\vrho}\,\psi^\vrho_\lambda(r)h_1(r).
    \]
    Similarly applying Equations \eqref{eq:Psi-der-zero} and \eqref{eq:Psi-topder}, we obtain
    \begin{equation}\label{eq:halfint-j1star}
      (f*\sh J(1))(0)=2\int_0^1dr\,r^{\vrho+1}(1-r)^{-1-\vrho}\Psi^\vrho(r)h_1(r)+\frac{(-1)^{\vrho-\frac12}2^{3-\vrho}\pi}{\Gamma\Parens1{\frac12-\vrho}}h_2(0).
    \end{equation}
    As $f(k\cdot 0)=f(0)$, Equation \eqref{eq:log-int-2} implies that 
    \[
      \frac{(-1)^{\vrho-\frac12}2^{3-\vrho}\pi}{\Gamma\Parens1{\frac12-\vrho}}h_2(0)=-2^{2(1-\vrho)}\pi f(0).
    \]

    To prove Equation \eqref{eq:fouinv-origin} and thereby complete the theorem's proof, it is therefore sufficient to see that after inserting the definition of $\Psi^\vrho(r)$, the order of integration over $r$ and $\lambda$ may be interchanged in the first summand on the right-hand side of Equation \eqref{eq:halfint-j1star}. But Equation \eqref{eq:sphfn-wt-der-est} implies the existence of a constant $C>0$ \scth
    \[
      \Abs1{\psi_\lambda^\vrho(r)}\sle C(1-r)^{-\frac{\vrho}2}
    \]
    \fa $\lambda\in i\ger a_\reals^*$ and $r\in[0,1)$. This gives sufficient bounds at $r=1$. On the other hand,  Equations \eqref{eq:psi-der} and \eqref{eq:psi-der-zero}, and the proof of the latter, show that $r^{\vrho+1}\psi_\lambda^\vrho(r)$ is uniformly bounded in $\lambda$ in a neighbourhood of $r=0$. Thus, the Fubini theorem applies, finally completing the theorem's proof. 
  \end{proof}

  \begin{Rem}
    The strategy of proof for \thmref{Th}{fouinv} in the case of $\vrho<0$ half-integral generalizes ideas from  Zirnbauer's proof of the corresponding result \cite[Theorem 2 and Corollary]{zirn-cmp} for the super-Poincar\'e disc, where $\vrho=-\tfrac12$. However, in order to implement the strategy, we needed sharper estimates on the spherical functions and knowledge of their derivatives. Note also that the determination of the constant depends on \thmref{Prop}{norm-const}, which required the application of localization techniques in its proof. 
  \end{Rem}


\end{document}